\documentclass[11pt]{article}
\usepackage{pgf}
\usepackage{lipsum}
\usepackage{amsfonts}
\usepackage{graphicx}
\usepackage{epstopdf}
\usepackage{algorithmic}
\ifpdf
  \DeclareGraphicsExtensions{.eps,.pdf,.png,.jpg}
\else
  \DeclareGraphicsExtensions{.eps}
\fi
\usepackage{color}
\usepackage{amsmath}
\usepackage{amsthm}
\usepackage{amssymb}
\usepackage{epsfig}
\usepackage{psfrag}
\usepackage{graphicx}
\usepackage{textcomp}
\usepackage{pgf}   
\usepackage{upgreek} 
\usepackage{mathrsfs,cite}

\newtheorem{lemma}{Lemma}
\newtheorem{proposition}{Proposition}

\newtheorem{remark}{Remark}
\newtheorem{example}{Example}
\newtheorem{corollary}{Corollary}

\newcommand{\R}{\mathbb{R}}
\newcommand{\N}{\mathbb{N}}

\newcommand{\lineunder}[2]{\LU{\begin{array}[t]{c}\underbrace{#1}\vspace*{.5em}\end{array}}{\mbox{\footnotesize\rm #2}}}

\newcommand{\LU}[2]{\begin{array}[t]{c}#1\vspace*{-1em}\\_{#2}\end{array}}

\renewcommand{\d}{{\rm d}}

\newcommand\Frakd{{\large\mbox{$\mathfrak{d}$}}}

\newcommand{\tto}{\buildrel {_{_{_{\mbox{$\hspace*{.5mm}\rightarrow$}}}}}\over{_{\mbox{$\rightarrow$}}}}
\newcommand\rca{{\rm rca}}
\newcommand\rba{{\rm rba}}
\newcommand\prca{{\rm rca}_1^+}
\newcommand\prba{{\rm rba}_1^+}
\newcommand\NNU{\boldsymbol\nu}
\newcommand\MMU{\boldsymbol\mu}
\newcommand\CHI{\boldsymbol\chi}
\newcommand\DELTA{\boldsymbol\delta}
\newcommand{\barOmega}{\mathchoice{\hspace*{.12em}{\overline{\hspace*{-.12em}\varOmega}}}
                                  {\hspace*{.12em}{\overline{\hspace*{-.12em}\varOmega}}}
                   {{\footnotesize\hspace*{.12em}{\overline{\hspace*{-.12em}\varOmega}}}}
                  {{\footnotesize\hspace*{.15em}{\overline{\hspace*{-.15em}\varOmega}}}}}

\newcommand\TTT{\color{black}}
\newcommand\EEE{\color{black}}

\usepackage{amsopn}

\begin{document}
\begin{sloppypar}

  \begin{center}
    {\LARGE\bf
      Fine metrizable convex relaxations of\\[.5em]
  parabolic optimal control problems} \footnote{The author acknowledges a partial support of the
CSF (Czech Science Foundation) project
19-04956S,
the M\v SMT \v CR (Ministery of Education of the Czech Rep.) project
CZ.02.1.01/0.0/0.0/15-003/0000493,
and the institutional support RVO: 61388998 (\v CR).}
  \end{center}
\bigskip\bigskip

\centerline{\large\sc Tom\'a\v s Roub\'\i\v cek}

\bigskip

\begin{center}
  {Mathematical Institute, Charles University,\\
Sokolovsk\'a 83, CZ-186~75~Praha~8, Czech Republic,
\\
and\\
Institute of Thermomechanics of the Czech Academy of Sciences,\\ 
Dolej\v skova 5, CZ-182 00 Praha 8, Czech Republic}

\end{center}

  

\bigskip

\begin{minipage}[t]{15cm}{\small
  \noindent{\bf Abstract.}
  Nonconvex optimal-control problems governed by evolution
problems in infinite-dimensional spaces (as e.g.\
parabolic boundary-value problems) needs a continuous
(and possibly also smooth) extension on some (preferably convex)
compactification, called relaxation, to guarantee existence
of their solutions and to facilitate analysis by relatively
conventional tools. When the control is valued in some subsets
of Lebesgue spaces, 
the usual extensions are either too coarse (allowing in fact
only very restricted nonlinearities) or too fine (being
nonmetrizable). To overcome these drawbacks,
a compromising convex compactification is here devised,
combining classical techniques for Young measures with Choquet theory.
This is applied to parabolic optimal control problems as far as
existence and optimality conditions concerns.

\medskip


\medskip

  }
\end{minipage}

\def\bfu{\boldsymbol u}
\def\bfr{\boldsymbol r}
\def\bfy{\boldsymbol y}
\def\bff{\boldsymbol f}
\def\bfs{\boldsymbol s}
\def\bfh{\boldsymbol h}
\def\bfv{\boldsymbol v}
\def\bfb{\boldsymbol b}
\def\bfq{\boldsymbol q}
\def\bfz{\boldsymbol z}
\def\BFS{\boldsymbol S}
\def\bfS{S}
\def\bfSp{S_p}
\def\bfU{\boldsymbol U}
\def\bfQ{\boldsymbol Q}
\def\bfpi{\boldsymbol\pi}
\def\LCS{L^1(\varOmega;C(B))}
\def\calR{\mathscr R}
\def\calY{\boldsymbol{\mathcal Y}}
\def\clR{\mathcal R}
\def\overlinebff{\hspace*{.15em}\overline{\hspace*{-.15em}\bff}}
\def\overlinebfh{\hspace*{.1em}\overline{\hspace*{-.1em}\bfh\hspace*{-.05em}}\hspace*{.05em}}
\def\overlinebfu{\hspace*{.1em}\overline{\hspace*{-.1em}\bfu\hspace*{-.05em}}\hspace*{.05em}}
\def\overlinebfpi{\hspace*{.15em}\overline{\hspace*{-.15em}\bfpi}}
\def\overlinebfQ{\hspace*{.15em}\overline{\hspace*{-.15em}\bfQ}}

\bigskip

\section{Introduction}\label{sect-int}

Relaxations in optimal control theory means usually a certain natural
extension of optimization problems. The adjective ``natural'' means
most often ``by continuity''.  The essential attribute of the
extended (called relaxed) problems is compactness of the set of admissible
relaxed controls, which ensures
existence and stability of solutions. An additional attribute is
convexity of this set, which allows for further analysis leading
to optimality conditions. A general theory of so-called convex
compactifications can be found in \cite{Roub20ROTV}.

A particular situation, which this paper is focused on, appears in
optimal control of evolution problems in infinite-dimensional spaces.
This abstract situation covers in particular optimal control of
systems governed by parabolic partial differential equations. 

Conventional relaxation in control theory of such evolution problems
uses the original controls ranging over an abstract topological space $S$ and
works with
continuous nonlinearities. After relaxation, this gives rise to a standard
$\sigma$-additive functions (measures) on the Borel $\sigma$-algebra
of Borel subsets. These measures are parameterized by time and possibly,
in the parabolic systems, also by space.
Such parameterized measures are called Young measures \cite{young1},
although L.C.\,Young worked rather with functionals because the measure
theory was rather only developing. The spirit of Young measures
as functionals (allowing more straightforwardly for various generalizations
or approximation) is accented in \cite{Roub20ROTV}.

In abstract evolution problems, the set $S$ where the
controls are valued is a compact metric space or, a bit more generally,
a so-called Polish space (separable completely metrizable topological space)
which is compact. This leads to Young measures parameterized by time
valued in conventional (i.e.\ $\sigma$-additive) probability measures supported 
on $S$, cf.\ e.g.\ \cite{ahmed,ahmed-2,avgerinos-papageorgiou,AvePaP90OCRC,CaRFVa04YMTS,EdmThi05ROCP,papageorgiou,XiaAhm93PRTE,Warg72OCDF}.
A modification for metrizable locally compact \TTT sets \EEE
was devised by \cite{Bald00NFYM}.

Often, $S$ is a subset of an infinite dimensional Banach space. 
In parabolic problems
interpreted as evolution problems on Sobolev spaces over
a domain $\varOmega\subset\R^d$, $d\in\N$, the set $S$ where the
controls are valued is typically a subset
of some Lebesgue space over $\varOmega$. More specifically, let us
consider
\begin{align}\label{S}
\bfSp=\{u\in L^p(\varOmega;\R^m);\ \ u(x)\in B
\ \text{ for a.a. }\ x\in\varOmega\}\,,
\end{align}
which is compact in its weak (or weak*) topology if $B\subset\R^m$ is bounded
and closed; \TTT except Remark~\ref{rem-DM}, the set $B$ will always be
bounded and $p$ only denotes some (or equally any) number such that
$1\le p<+\infty$ and wants to emphasize that (except  Remark~\ref{rem-DM})
this subset of $L^\infty(\varOmega;\R^m)$ is considered in the
$L^p$-topology with $p\ne\infty$. \EEE
The nonlinearities occurring in concrete optimal-control problems have typically
a local form of the type $u\mapsto h(x,u(x))$ with some Carath\'eodory
mapping $h:\varOmega\times\R^m\to\R^k$. Yet, in this space, such mappings 
are weakly continuous on $\bfSp$ only if $h(x,\cdot)$ is affine for a.a.\
$x\in\varOmega$. This hidden effect makes the approach from
\cite{avgerinos-papageorgiou,AvePaP90OCRC,papageorgiou,Warg72OCDF,XiaAhm93PRTE}
in fact very restrictive, as 
functions of controls \TTT which are not affine \EEE
do not admit continuous extension in terms of Young measures.

A finer convex compactification was devised by Fattorini
\cite{fattorini,Fatt94ETMP,fattorini-4}, allowing for a general continuous
nonlinearities on
$\bfSp$ but using the rather abstract concept of Young measures valued in
probability regular finite additive measures ``rba'' on $\bfSp$, or equivalently
\cite{Roub20ROTV} as standard probability regular countably additive
measures ``rca'' but on the \v Cech-Stone compactification $\beta\bfSp$ of
$\bfSp$. Such compactification is not metrizable, and one cannot work with
conventional sequences but, instead, the general-topological concept
of nets and Moore-Smith convergence must be used.

The goal of this paper is to devised a compromising relaxation which
admits a wider class of nonlinearities than only affine while still working
with conventional $\sigma$-additive measures and conventional sequences.
For this, a characterization
of extreme Young measures together with celebrated Choquet-Bishop-de\,Leeuw
\cite{BisDeL59RLFM,Choq60TRIE} theory is used first in the ``static''
situation in Section~\ref{sec-YM}, and then for the evolution situations
parameterized by time in Section~\ref{sec-YM-time}. Eventually,
in Section~\ref{sec-appl}, application to optimal control of parabolic
partial differential equations is briefly shown.

 \section{Young measures and probability
measures on them}\label{sec-YM}

 Let us begin with some definitions and brief presentation of basic needed
 concepts and facts. An algebra on a set $M$ is a collection of subsets
of $M$ closed on the complements and finite unions, including also an
empty set. If it is also closed on union of countable number of sets,
then it is called an $\sigma$-algebra.
We denote by $C(M)$ a space of continuous {\it bounded} function on a
topological space $M$. In fact, it is algebraically also an algebra and,
if $M$ is compact,
it is a Banach space and, by the classical Riesz theorem, its dual $C(M)^*$ is
isometrically isomorphic to the Banach space of Borel measures denoted by
$\rca(M)$, i.e.\ of regular bounded countably-additive (so-called
$\sigma$-additive) set functions on the $\sigma$-algebra of Borel subsets of
$M$. 

If $M$ is a (not necessarily compact) normal topological space, $C(M)^*$ is
isometrically isomorphic to the Banach space of regular bounded finitely
additive (not necessarily $\sigma$-additive) set functions on the
algebra generated by all subsets of $M$, denoted by $\rba(M)$. Actually,
$\rba(M)=\rca(\beta M)$ with $\beta M$ the \v Cech-Stone compactification
$\beta M$ of $M$.

The subsets of $\rca(M)$ and $\rba(M)$ consisting from positive measures
having a unit mass (i.e.\ probability measures) will be denoted by
$\prca$ and $\prba(M)$, respectively.
If $C(M)$ is separable, then the weak*
topology on $\prca(M)$ or $\prba(M)$ is metrizable.

If $M$ is a domain in an Euclidean space equipped
with the Lebesgue measure (denoted then mostly by $\varOmega\subset\R^d$),
then $L^p(M;\R^n)$ will denote the Lebesgue space of all
measurable $\R^n$-valued whose $p$-power is integrable.

An important attribute of $\prca(M)$ and of $\prba(M)$
is convexity.
Let us remind that a point $z$ in a convex set $K$ is
called {\it extreme} in $K$ if there is no open interval
in $K$ containing $z$; in other words, $z=az_1+(1{-}a)z_2$ for some
$a\in(0,1)$ and some $z_1,z_2\in K$ implies $z_1=z_2$. 
\TTT The set of the extreme points will be denoted as ${\rm ext}\,K$. \EEE
Let us note that the set of all extreme points of a
metrizable convex compact $(K,\TTT\rho\EEE)$ is a Borel set (more precisely a
$G_\delta$-set), being a countable intersection
of open sets as complements to the closed set $\{z\in K;\ \exists z_1,z_2\in K:\
z=\frac12z_1+\frac12z_2\ \&\ \TTT\rho\EEE(z_1,z_2)\le\epsilon\}$
for $\epsilon>0$.

One of the important ingredients used below is that every point $z$ of  a
convex compact set $K$ is an average of the extreme points according
to a certain probability measure $\mu$ supported on extreme points in the
sense
\begin{align}
\forall z\in K\ \ \exists\,\mu\in\prca({\rm ext}\,K)\ \
\forall\,f\in{\rm Aff}\,K:\quad
f\Big(\int_{{\rm ext}\,K}\!\!\!\!\widetilde{z}\,\mu(\d \widetilde{z})\Big)
=\int_{{\rm ext}\,K}\!\!\!\!f(\widetilde{z})\,\mu(\d \widetilde{z})=f(z)\,,
\label{ext+}\end{align}
where ${\rm Aff}\,K$ denotes the set of all affine continuous functions on $A$.
In other words, any $z\in K$ is a so-called barycentre of a probability
measure supported on ${\rm ext}\,K$. This is known as a
{\it Choquet-Bishop-de\,Leeuw representation theorem}\label{Choquet-theorem}
\cite{BisDeL59RLFM,Choq60TRIE}; cf.\  also
e.g.\ \cite{Alfs71CCSB,LMNS10IRT}.
Recall that $z\in K$ is called a barycentre of
$\mu\in\prca(K)$ if $f(z)=\int_Kf(\widetilde{z})\,\mu(\d \widetilde{z})$
for any affine continuous $f:K\to\R$, so that the last equation in
\eqref{ext+} says that $z$ is a barycentre of $\mu$.

We now briefly recall the {\it classical Young measures}.
We consider the set $\bfSp$ from \eqref{S}
with $B\subset\R^m$ compact, not necessarily convex. As $B$ is bounded, the set
$\bfSp$ actually does not depend on $1\le p\le+\infty$. When endowed with the
norm topology from $L^p(\varOmega;\R^m)$,
it becomes a normal topological space.
This topology is separable and does not depend on $1\le p<\infty$, but
$\bfS_\infty$ has a strictly finer (and non-separable) topology.

The notation $L^\infty_{\rm w*}(\varOmega;X^*)$ stands for the Banach space of
weakly* measurable mappings $\nu:\varOmega\to X^*$ for some Banach space $X$,
i.e.\ $x\mapsto\langle\nu(x),h(x)\rangle$
is measurable
for any $h\in L^1(\varOmega;X)$. Here we use it for $X=C(B)$ and later
also for some subspaces of $C(\bfSp)$. By the
Dunford-Pettis' theorem combined with the mentioned Riesz theorem,
$L^1(\varOmega;C(B))^*\cong L^\infty_{\rm w*}(\varOmega;\rca(B))$.
For $\nu\in L^\infty_{\rm w*}(\varOmega;\rca(B))$, it is customary to write
$\nu_x$ instead of $\nu(x)$. We define the set of Young measures 
\begin{align}\label{Young}
{\mathcal Y}(\varOmega;B):=\big\{\nu\in L^\infty_{\rm w*}(\varOmega;\rca(B));
\ \nu_x\in\prca(B)\ \text{ for a.a. }\,x\in\varOmega\big\}\,.
\end{align}
It is obvious that ${\mathcal Y}(\varOmega;B)$ is convex,
weakly* compact, and metrizable. The set $\bfSp$ is embedded into
${\mathcal Y}(\varOmega;B)$ by the mapping $[\delta(u)]_x=\delta_{u(x)}$
where $\delta_s\in\prca(B)$ denotes the Dirac measure supported at $s\in B$.
By a direct construction of fast oscillating sequences, one can show that
this embedding is weakly* dense and thus, in particular,
${\mathcal Y}(\varOmega;B)$ is separable.
It is important that
the embedding $\delta:\bfSp\to{\mathcal Y}(\varOmega;B)$ is even
(strong,weak*)-homeomorphical with respect to the strong topology of
$L^p(\varOmega;\R^m)$ for any $1\le p<+\infty$, although not for $p=+\infty$.
Here we note that $\delta(u_k)\to\delta(u)$ weakly* in 
${\mathcal Y}(\varOmega;B)$ implies, when tested by
\begin{align}\label{test}
  \big(h:(x,z)\mapsto|z{-}u(x)|^p\big)\in L^1(\varOmega;C(B))\,,
\end{align}
that
$\langle\delta(u_k){-}\delta(u),h\rangle=\int_\varOmega|u_k{-}u|^p\d x\to0$.
\TTT Let us remind that $B$ is considered bounded (and closed) and
$p<+\infty$, otherwise the inclusion in \eqref{test} would not hold. \EEE

The other important ingredient used below is that each extreme point
$\nu=\{\nu_x\}_{x\in\varOmega}$ in the set of all Young measures
${\mathcal Y}(\varOmega;B)$ \TTT is \EEE composed from
Diracs, i.e.\ $\nu_x=\delta_{u(x)}$ for a.a.\ $x\in\varOmega$ with some
$u\in\bfSp$; \TTT see \EEE 
Berliocchi and Lasry \cite[Proposition~II.3]{berliocchi-lasry} or Castaing and
Valadier \cite[Thm.\,IV.15]{castaing-valadier}, cf.\ also \cite{kruzik-tr}.
The extreme points of ${\mathcal Y}(\varOmega;B)$ are thus a dense and,
as mentioned above, $G_\delta$-set in ${\mathcal Y}(\varOmega;B)$.

\def\MU{\mu}

\begin{lemma}\label{lem1}
\TTT Let $B\subset\R^m$ is compact and $1\le p<+\infty$. \EEE 
Any Young measure $\nu\in{\mathcal Y}(\varOmega;B)$ can be represented
(in a non-unique way in general)
by a probability measure $\MU
$ supported on $\bfSp$. More specifically,
\begin{align}\nonumber
\forall\nu\in{\mathcal Y}(\varOmega;B)\ \ \ 
&\exists\,\MU\in\prca({\mathcal Y}(\varOmega;B)),\ \
{\rm supp}\,\MU\subset\Frakd(\bfSp)\ \ \ \forall
h\in\LCS:
\\&\qquad\qquad
\int_\varOmega\!\int_{B}h(x,z)\,\nu_x(\d z)\,\d x
=\int_{\bfSp}\!\int_\varOmega h(x,u(x))\,\d x\,\MU(\d u)\,,
\label{ext++}
\end{align}
where we identified $\MU(\Frakd(A))$ and $\MU(A)$ for
$A\subset\bfSp$; here
$\Frakd:\bfSp\to\prca(\bfSp):u
\mapsto\Frakd_u^{}$ with $\Frakd_u^{}\in\prca(\bfSp)$ denoting
the Dirac measure supported at $u\in \bfSp$.
Thus, in fact, $\MU\in\prca(\bfSp)$.
\end{lemma}

This means that $\Frakd_u^{}$ as a functional on $C(\bfSp)$ defined
by $v\mapsto v(u)$ for any $v\in C(\bfSp)$, in contrast to
$\delta(u)=\{\delta_{u(x)}\}_{x\in\varOmega}^{}\in{\mathcal Y}(\varOmega;B)$.

\begin{proof}[Proof of Lemma~\ref{lem1}]
In view of the abstract result \eqref{ext+}, i.e.\ the Choquet-Bishop-de\,Leeuw
representation theorem applied on the convex compact $K={\mathcal Y}(\varOmega;B)$.
We thus obtain a probability measure
$\MU\in\prca({\mathcal Y}(\varOmega;B))$
supported on ${\rm ext}({\mathcal Y}(\varOmega;B))$.
As mentioned above, ${\rm ext}({\mathcal Y}(\varOmega;B))=\Frakd(\bfSp)$.
Since ${\mathcal Y}(\varOmega;B)$ is metrizable, $\Frakd(\bfSp)$
is a Borel subset in ${\mathcal Y}(\varOmega;B)$, and $\MU$
is a Borel measure on it. Also realize that any weakly*
continuous affine function on $K={\mathcal Y}(\varOmega;B)\subset
L^\infty_{\rm w*}(\varOmega;\rca(B))\cong L^1(\varOmega;C(B))^*$ is of
the form $\nu\mapsto\int_\varOmega\int_{B}h(x,z)\,\nu_x(\d z)\,\d x$ for some
$h\in L^1(\varOmega;C(B))$. Thus \eqref{ext+} yields \eqref{ext++}.

It is important that, as mentioned above, the embedding $\Frakd:\bfSp$
is (strong,weak*)-homeomorphical so that
the weak* topology on $\Frakd(\bfSp)$ induces just the strong topology
on $\bfSp$. Thus the measure on $\Frakd(\bfSp)$ induces a Borel measure
on $\bfSp$, referring to the Borel $\sigma$-algebra on $\bfSp$ with respect to
the $L^p$-norm, $p<+\infty$, such a measure on $\bfSp$ being again denoted by
$\MU$.
\end{proof}

Let us still remind a canonical construction of compactifications, here
applied to $\bfSp$. For this, we consider a general complete closed sub-ring
$\calR$ of $C(\bfSp)$ containing constants.
Every such a ring $\calR$ is also a commutative Banach algebra and determines
a compactification $\gamma_\calR^{}\bfSp$ of $\bfSp$ as a subset of
$\calR^*$ endowed with the weak* topology consisting of multiplicative
means, i.e.\
$$
\gamma_\calR^{}\bfSp:=\big\{\mu\in\calR^*;\ \|\mu\|=1,\ \mu(1)=1,\ 
\forall v_1,v_2\in\calR:\ \langle\mu,v_{1}v_2\rangle
=\langle\mu,v_1\rangle\,\langle\mu,v_2\rangle\big\}\,.
$$
This means $\gamma_\calR^{}\bfSp$ is compact when endowed by the weak* topology
of $\calR^*$ and the embedding $e:\bfSp\to\gamma_\calR^{}\bfSp$ defined
by $\langle e(u),v\rangle=v(u)$ is homeomorphical.
Let us recall that a ring is called
complete if it separates closed subsets of $\bfSp$ from points in $\bfSp$.
This means that, for any
$A\subset\bfSp$ closed and $u_0\not\in\bfSp\setminus A$, there is
$v\in C({\mathcal Y}(\varOmega;B))|_{\bfSp}^{}$ such that $v(u_0)=0$ and
$v(A)=1$.
The functionals from $\gamma_\calR^{}\bfSp$ are positive in the sense
that $\mu(v)\ge0$ for any $v\in\calR$ with $f(\cdot)\ge0$ on $\bfS$.
Each
function from $\calR$ admits a (uniquely determined) continuous extension on
$\gamma_\calR^{}\bfSp$. Thus $\calR$ is isometrically isomorphic with the space
$C(\gamma_\calR^{}\bfSp)$. By the Riesz theorem,
$\calR^*\cong\rca(\gamma_\calR^{}\bfSp)$ and
\begin{align}\label{duality}
\langle\mu,v\rangle=\int_{\gamma_\calR^{}\bfSp}\!\!\!\!\!\!\overline v(s)\mu(\d s)
\ \ \ \text{ with }\,\overline v\in C(\gamma_\calR^{}\bfSp)\,\text{ a continuous
  extension of }\,v\in\calR\,.
\end{align}

As already mentioned in Sect.\,\ref{sect-int}, the construction
$\prba(\bfSp)\cong \prca(\beta \bfSp)$ is non-metrizable (and thus rather
constructive) because it is based on the nonseparable space of test functions
$C(\bfSp)$. It is thus desirable to consider some subspace of $C(\bfSp)$ which
would be separable but still bigger than ${\rm Aff}(\bfSp)$. Motivated by
Lemma~\ref{lem1}, we take the choice
\begin{align}\label{tests}
  \calR=
C({\mathcal Y}(\varOmega;B))\big|_{\bfSp}\!=
\big\{v\in C(\bfSp);\ \exists\,\overline v\in C({\mathcal Y}(\varOmega;B)):\
 v=\overline v\circ\delta\big\}\,.
\end{align}
This is obviously a sub-ring of $C(\bfSp)$ containing constants.

\begin{lemma}\label{lem2}
The ring $\calR$ from \eqref{tests} is complete and separable and the
compactification $\gamma_\calR^{}\bfSp$ is metrizable and homeomorphical with
${\mathcal Y}(\varOmega;B)$.
\end{lemma}

\begin{proof}
Since ${\rm dist}(u,A)=\epsilon>0$, one can take 
$v(u)=\min(\epsilon,\|u{-}u_0\|_{L^p(\varOmega;\R^m)}^p)$
which can indeed be continuously extended on ${\cal Y}(\varOmega;B)$
as $v(\nu)=\min(\epsilon,\langle\nu,h\rangle)$ with
$h$ from \eqref{test}; note that such $h$ is an integrand from
$L^1(\varOmega;C(B))$.
In fact, $C({\mathcal Y}(\varOmega;B))|_{\bfSp}^{}$ is the smallest closed ring
containing $\varPsi(L^1(\varOmega;C(B)))$,
 where $\varPsi$ is a linear
operator from $L^1(\varOmega;C(B))$ to $C(\bfSp)$ defined by
\begin{align}\label{Psi}
\varPsi h:u\mapsto\int_\varOmega h(x,u(x))\,\d x\,.
\end{align}

Let us also remind that ${\mathcal Y}(\varOmega;B)$ is a metrizable separable
compact.
Hence $C({\mathcal Y}(\varOmega;B))$ itself
is separable. The separability holds also for \eqref{tests}.

Since $\calR$ is separable, bounded sets in 
its dual endowed with the weak* topology 
(and in particular $\gamma_\calR^{}\bfSp$) are metrizable.
The homeomorphism between $\gamma_\calR^{}\bfSp$ and
${\mathcal Y}(\varOmega;B)$ is realized by the adjoint 
operator to the embedding of $\varPsi(L^1(\varOmega;C(B)))$ into
$C({\mathcal Y}(\varOmega;B))|_{\bfSp}^{}$;
here it is important that $\varPsi(L^1(\varOmega;C(B)))$ is a
so-called convexifying subspace of $C(\bfSp)$ in the sense
that any $u_1,u_2\in\bfSp$ admits a sequence $\{u_k\}_{k\in\N}$ such that
$f(u_1){+}f(u_2)=2\lim_{k\to\infty}f(u_k)$ for any 
$f\in\varPsi(L^1(\varOmega;C(B)))$, cf.\ \cite[Sect.\,2.2 and 3.1]{Roub20ROTV}. 
\end{proof}

\begin{example}\upshape
Let us still illustrate Lemma~\ref{lem1} on a piece-wise
homogeneous two-atomic Young measure
\begin{align}\label{example-NU}
\nu_x=\begin{cases}
\frac12\delta_{u_1^{}(x)}+\frac12\delta_{u_2^{}(x)}&\text{ for }x\in A\,,\\
\frac14\delta_{u_1^{}(x)}+\frac34\delta_{u_2^{}(x)}&\text{ for }x\in \varOmega{\setminus}A\,,
\end{cases}
\end{align}
with some $u_1\ne u_2$ and $A\subset\varOmega$ measurable. Then $\MU$ from
Lemma~\ref{lem1} takes (for example) the form
\begin{align}\label{example-NNU}
  \ \ \ \ \MU=a\Frakd_{u_{11}^{}}^{}\!\!
+\Big(\frac12{-}a\Big)\Frakd_{u_{12}^{}}^{}\!\!
+\Big(\frac14{-}a\Big)\Frakd_{u_{21}^{}}^{}\!\!
+\Big(\frac14{+}a\Big)\Frakd_{u_{22}^{}}^{}
\end{align}
with an arbitrary parameter $0\le a\le1/4$ and with $u_{11}=u_1$, $u_{22}=u_2$,
\begin{align*}
u_{12}(x)=\begin{cases}u_1(x),&\\[-.3em]
u_2(x),&
\end{cases}\ \ \ \text{ and }\ \ \ \
u_{21}(x)=\begin{cases}u_2(x)&\text{ for }x\in A\,,\\[-.3em]
u_1(x)&\text{ for }x\in \varOmega{\setminus}A\,.
\end{cases}
\end{align*}
For $a=0$ and for $a=1/4$, the four-atomic measure \eqref{example-NNU}
degenerates
to only three-atomic ones. In particular, it illustrates non-uniqueness of the
probability measure from Lemma~\ref{lem1}. Even more,
\eqref{example-NNU} does not cover all representations of $\nu$ from
\eqref{example-NU}.
Although these measures cannot be distinguished when tested by test functions
from $\varPsi(L^1(\varOmega;C(B)))$, they can be distinguished from each other
when tested by functions from $C({\mathcal Y}(\varOmega;B))|_{\bfSp}^{}$; for
example, if $\varrho$ is a metric of ${\mathcal Y}(\varOmega;B)$, one can
take $\varrho(\cdot,\MU)|_{\bfSp}^{}$ with $\MU$ from
\eqref{example-NNU} for some specific  $0\le a\le1/4$.
\end{example}

\begin{remark}[{\sl An approximation of Young measures}]\label{rem-approx}
  \upshape
Various numerical schemes have been devised to numerical approximation
of Young measures, cf.\ \cite{Roub20ROTV} for a survey. Lemma~\ref{lem1}
inspires an approximation by a convex combination of elements from $\bfSp$.
Actually, this sort of approximation is supported by arguments that
each element of convex compact sets (i.e.\ here the set of Young measures)
can be approximated by a convex combination of extreme points due to the 
celebrated Kre\u{\i}n-Milman theorem.
Here one can consider a fixed countable
collection $\{u_l\}_{l\in\N}^{}$ dense in $\bfSp$ and, for any $\ell\in\N$, define
the finite-dimensional convex subset of $\prca(\bfSp)$ as
\begin{align}
\Big\{\MU=\sum_{l=1}^\ell a_l\Frakd_{u_l}^{};\
\exists\{a_l\}_{l=1}^{\ell}\,,\ \ a_l\ge0\,,\ \ 
&\sum_{l=1}^{\ell}\!a_l=1\,\Big\}\,.
\label{(3.5-mix)}\end{align}
This approximation, devised by V.M.\,Tikhomirov \cite{Tikh82LPOC}, was used
e.g.\ in \cite{AvaMag14MCPM,AvaMag17RCOC} under the name a
{\it mix of controls}. Passing
$\ell\to\infty$, the sets \eqref{(3.5-mix)} increase and their union is dense
in $\prca(\bfSp)$ due to the weak* density of $\{\Frakd_{u_l}\}_{l\in\N}^{}$.
This allows
for the convergence proof behind this sort of convex approximation.
\end{remark}

\begin{remark}[{\sl Special probability measures on $\bfSp$}]
\upshape
We do not claim that each $\MU\in\prca(\bfSp)$ corresponds
to some $\nu\in{\mathcal Y}(\varOmega;B)$ via \eqref{ext++}. For further
purposes, let us denote the set of such ``special'' $\MU$'s by
\begin{align}\nonumber
{\rm srca}_1^+(\bfSp):=\bigg\{\MU\in\prca(\bfSp);&\ \ 
  \exists \nu\in{\mathcal Y}(\varOmega;B)\ \ \forall
  h\in\LCS: 
  \\[-.7em]&
\int_{\bfSp}\![\varPsi h](u)\,\MU(\d u)=
  \int_\varOmega\!\int_{B}h(x,z)\,\nu_x(\d z)\,\d x  \,\bigg\}\,.
\label{ext+++}
\end{align}
Although this set is defined only very implicitly, we can nevertheless
see that the set ${\rm srca}_1^+(\bfSp)$ is convex. Indeed,
for $\MU_1,\MU_2\in\prca(\bfSp)$, we can take
$\nu_1,\nu_2\in{\mathcal Y}(\varOmega;B)$ such that, for $i=1,2$, it holds
\begin{align}
  \!\!\forall h\!\in\!\LCS:\ \ 
  \int_{\bfSp}\![\varPsi h](u)\,\MU_i(\d u)=
  \int_\varOmega\!\int_{B}h(x,z)\,\big[\nu_i\big]_x^{}(\d z)\,\d x
  \,.
\label{srca+}\end{align}
As ${\mathcal Y}(\varOmega;B)$ is convex, also
$\nu=\frac12\nu_1{+}\frac12\nu_2\in{\mathcal Y}(\varOmega;B)$.
Thus $\MU=\frac12\MU_1{+}\frac12\MU_2$ satisfies the identity in \eqref{ext+++}.
\end{remark}

\begin{example}[{\sl
Special functions from $C({\mathcal Y}(\varOmega;B))|_{\bfSp}^{}$}]
\label{rem-1}\upshape
Using $\varPsi$ from \eqref{Psi}, let us consider
\begin{align}\label{tests+}
  \calR_0:=\Big\{v{\in}C(\bfSp);\ \exists\,
m,n{\in}\N,\ f_{ij}{\in} C(\R),\ h_{ij}{\in} L^1(\varOmega;C(B)){:}\
  v=\sum_{i=1}^m\prod_{j=1}^nf_{ij}{\circ}\varPsi h_{ij}\Big\}.
\end{align}
This is a complete separable ring. Each $v\in\calR_0$ admits a weakly*
continuous extension $\overline{v}$ on ${\mathcal Y}(\varOmega;B)$, given
explicitly by
$$
\overline{v}(\nu)=
\sum_{i=1}^m\prod_{j=1}^nf_{ij}
\Big(\int_\varOmega\int_Bh_{ij}(x,z)\,\nu_x(\d z)\d x\Big)
$$
for all $\nu\in{\mathcal Y}(\varOmega;B)$. Thus $\calR_0\subset
C({\mathcal Y}(\varOmega;B))|_{\bfSp}^{}$. Moreover, by Lemma~\ref{lem1}, there
is $\mu\in\prca(\bfSp)$ depending on $\nu$ such that
$$
\overline{v}(\nu)=\sum_{i=1}^m\prod_{j=1}^nf_{ij}\Big(\int_{\bfSp}\!\!
\big[\varPsi h_{ij}\big](u)\,\mu(\d u)\Big)
\,.
$$
Since $C({\mathcal Y}(\varOmega;B))|_{\bfSp}^{}$ is the smallest closed ring
containing the linear space $\varPsi(L^1(\varOmega;C(B)))$,
$\calR_0$ is dense in $C({\mathcal Y}(\varOmega;B))|_{\bfSp}^{}$.
\end{example}

\def\LL{\boldsymbol L}
\def\HH{\boldsymbol H}

\section{Young measures parameterized by time}\label{sec-YM-time}

We will now extend the construction from Sect.\,\ref{sec-YM} to be
applicable for evolution problems. To this goal, we consider a time interval
$I=[0,T]$ with some fixed time horizon $T>0$.
We will use the standard notation  $L^p(I;X)$ for the Lebesgue-Bochner space of
abstract functions $I\to X$ whose $X$-norm is valued in $L^p(I)$.
From now on, let us agree to use boldface fonts for functions of time valued 
in spaces functions (or their duals) on $\varOmega$, or also functions of
such arguments. We will consider
the set of ``admissible controls''
\begin{align}\nonumber
\bfU_{\rm ad}&\,=\,\big\{u\in L^p(I{\times}\varOmega;\R^m);\
u(\cdot)\in B\text{ a.e.\ on }I{\times}\varOmega\big\}
\\&\,\cong\,
\big\{\bfu\in \LL^p(I;L^p(\varOmega;\R^m));\
\bfu(\cdot)\in\bfSp\text{ a.e.\ on }I\big\}
\label{Uad}\end{align}
with $\bfSp$ from 
\eqref{S} and with identifying
$\bfu(t,\cdot)$ and $\bfu(t)\in\bfSp$.

Let us start with a general construction, advancing the scheme devised
by H.\,Fattorini who used the non-separable space $\LL^1(I;C(\bfSp))$
of test functions. Here, instead of the whole (non-separable) ring
$C(\bfSp)$, we consider a general complete closed sub-ring $\calR$
of $C(\bfSp))$ containing constants, and then we can consider the test-function
space $\LL^1(I;\calR)$.
Instead of $\LL^1(I;C(\bfSp))$, we now suggest to use $\LL^1(I;\calR)$.
If $\calR$ is 
\TTT also \EEE separable, also $\LL^1(I;\calR)$ is separable and both
$\gamma_\calR^{}\bfSp$ and bounded sets in $\LL^1(I;\calR)^*$ are metrizable. 
Then, by Dunford-Pettis' theorem (as used in
\cite[Thms.\,12.2.4 and 12.2.11]{fattorini-4}), it holds
\begin{align}\label{(4.5.27+)}
\LL^1(I;\calR)^*\,\cong\,\LL^\infty_{\rm w*}(I;\calR^*)
\,\cong\, \LL^\infty_{\rm w*}(I;\rca(\gamma_\calR^{}\bfSp))\,.
\end{align}
Like before for $\nu(x)\equiv\nu_x$, we use the convention
$\NNU(t)\equiv\NNU_t$. The duality
between $\LL^\infty_{\rm w*}(I;\rca(\gamma_\calR^{}\bfSp))$ and
$\LL^1(I;\calR)$ is then
\begin{align}
\big\langle\NNU,\bfh\big\rangle
=\int_0^T\!\!\int_{\gamma_\calR^{}\bfSp}\!\!\!\overlinebfh(t,\bfs)\,\NNU_t(\d\bfs)\d t
\ \ \ \text{ for }\ \bfh\in \LL^1(I;
\calR)\,,
\label{embed+}\end{align}
where, like in \eqref{duality},
$\overlinebfh(t,\cdot)\in C(\gamma_\calR^{}\bfSp)$ is the (uniquely defined)
continuous extension of $\bfh(t,\cdot)\in\calR$ on $\gamma_\calR^{}\bfSp$.

Like in \eqref{Young}, we further define the set of Young measures with values
supported on $\gamma_\calR^{}\bfSp$ as
\begin{align}\label{Young+}
\calY(I;\gamma_\calR^{}\bfSp):=\{\NNU\in
\LL^\infty_{\rm w*}(I;\rca(\gamma_\calR^{}\bfSp));\
\NNU_t\in\prca(\gamma_\calR^{}\bfSp)\ \text{ for a.a. }t\in I\}\,.
\end{align}
Like \eqref{embed+}, the embedding $\DELTA:\bfU_{\rm ad}\to
\calY(I;\gamma_\calR^{}\bfSp)$
is defined as $\DELTA(\bfu)=\{\DELTA_{\bfu(t)}^{}\}_{t\in I}$, i.e.\ 
\begin{align}
\big\langle\DELTA(\bfu),\bfh\big\rangle
=\int_0^T\!\!\bfh(t,\bfu(t))\,\d t\ \ \ \text{ for }\
\bfh\in \LL^1(I;\calR)\,.
\label{embed-}\end{align}

\begin{proposition}\label{prop-0}
Let $\calR\subset C(\bfSp)$ be a separable complete closed sub-ring $\calR$
of $C(\bfSp)$ containing constants. Then the set of Young measures
$\calY(I;\gamma_\calR^{}\bfSp)$ from \eqref{Young+}
is convex, weakly* sequentially compact and separable.
Moreover, the embedding $\DELTA$ from \eqref{embed-} is
(strong,weak*)-continuous
if $\bfU_{\rm ad}$ is equipped with the strong topology
of $L^1(I{\times}\varOmega;\R^m)$ 
and
$\DELTA(\bfU_{\rm ad})$ is sequentially weakly* dense in
$\calY(I;\gamma_\calR^{}\bfSp)$.
\end{proposition}

\begin{proof}
The convexity and weak*-compactness of $\calY(I;\gamma_\calR^{}\bfSp)
\subset \LL^\infty_{\rm w*}(I;
\rca(\gamma_\calR^{}\bfSp))$ is obvious from the definition
 of the convex, closed, bounded set \eqref{Young+}. 
The continuity of $\DELTA$ follows from the continuity of the
Nemytski\u\i\ mappings induced by the integrands
$\bfh\in \LL^1(I;
\calR)$ in \eqref{embed-}.

The metrizability of the weak* topology on $\calY(I;\gamma_\calR^{}\bfSp)$
follows from the separability of $\LL^1(I;
\calR)$, relying on the
separability of $
\calR$.

The density of $\DELTA(\bfSp)$, i.e.\ 
the attainability of a general $\NNU\in\calY(I;\gamma_\calR^{}\bfSp)$
in the sense that $\NNU=\text{w*-}\lim_{k\to\infty}\DELTA(\bfu_k)$ for some
sequence $\{\bfu_k\}_{k\in\N}\subset \bfU_{\rm ad}$,
follows by the standard arguments used for conventional
Young measures, i.e.\ by an explicit construction of a sequence oscillating
fast in time, cf.\ \cite[Thm.\,12.6.7]{fattorini-4}
or \cite[Thm.\,3.6]{Roub20ROTV}. For this, an essential fact is that the
Lebesgue measure on $I$ is non-atomic. In particular, as $\bfSp$ is separable,
$\calY(I;\gamma_\calR^{}\bfSp)$  is separable, too.
\end{proof}

We now use this general construction for
the special choice $\calR=C({\mathcal Y}(\varOmega;B))|_{\bfSp}^{}$ as in
\eqref{tests}. We already 
\TTT showed \EEE in Lemma~\ref{lem2} that this ring is complete.
We thus consider the Banach space of
test functions $\LL^1(I;C({\mathcal Y}(\varOmega;B))|_{\bfSp}^{})$
and embed $\bfU_{\rm ad}$ into the dual of this test-function space
as in \eqref{(4.5.27+)} with $(C({\mathcal Y}(\varOmega;B))|_{\bfSp}^{})^*
\cong \rca({\mathcal Y}(\varOmega;B))$, i.e.\ into
\begin{align}\label{(4.5.27++)}
\big(\LL^1(I;C({\mathcal Y}(\varOmega;B))|_{\bfSp}^{}\big)^*
\,\cong\,
\LL^\infty_{\rm w*}(I;C({\mathcal Y}(\varOmega;B))|_{\bfSp}^*)
\,\cong\, \LL^\infty_{\rm w*}(I;\rca({\mathcal Y}(\varOmega;B)))\,.
\end{align}
By exploitation of Proposition~\ref{prop-0} with 
Lemma~\ref{lem2}, we have
\begin{align}
\calY(I;{\mathcal Y}(\varOmega;B))=\{&\NNU\in
\LL^\infty_{\rm w*}(I;\rca({\mathcal Y}(\varOmega;B)));
\ \ \NNU(t)\in\prca({\mathcal Y}(\varOmega;B))\,\text{ for a.a. }t{\in}I\}\,.
\label{special-YM}\end{align}

\begin{proposition}\label{prop-1}
The set of Young measures $\calY(I;{\mathcal Y}(\varOmega;B))$
is convex and a sepa\-rable metrizable compact, 
the embedding $\DELTA:\bfU_{\rm ad}\to\calY(I;{\mathcal Y}(\varOmega;B))$
again defined as \eqref{embed-} is (strong,weak*)-homeomorphical
if $\bfU_{\rm ad}$ is equipped with the strong topology
of $L^p(I{\times}\varOmega;\R^m)$ with any $1\le p<+\infty$
and $\DELTA(\bfU_{\rm ad})$ is weakly* dense in $\calY(I;{\mathcal Y}(\varOmega;B))$.
\end{proposition}

\begin{proof}
Most of the assertion follows from Proposition~\ref{prop-0} as a special
case for the choice $\calR=C({\mathcal Y}(\varOmega;B))|_{\bfSp}^{}$.

Considering a sequence $\{\bfu_k\}_{k\in\N}\subset\bfU_{\rm ad}$ and
$\bfu\in\bfU_{\rm ad}$ such that $\DELTA(\bfu_k)\to\DELTA(\bfu)$ weakly*,
we can prove that $\bfu_k\to\bfu$ strongly in $L^p(I{\times}\varOmega;\R^m)$.
Indeed, we can take $\bfh\in \LL^1(I;C({\mathcal Y}(\varOmega;B))|_{\bfSp}^{})$
defined as
$$
\bfh(t,\bfs)=\int_\varOmega\big|\bfs(x){-}[\bfu(t)](x)\big|^p\,\d x\,.
$$
This gives $\langle\DELTA(\bfu_k){-}\DELTA(\bfu),\bfh\rangle
=\|\bfu_k{-}\bfu\|_{L^p(I{\times}\varOmega;\R^m)}^p\to0$. Here it is important that
$\bfh(t,\bfs)=\int_\varOmega h(t,x,\bfs(x))\d x$ for 
$h(t,x,z)=|z{-}[\bfu(t)](x)|^p$ with $h(t,\cdot,\cdot)\in L^1(\varOmega;C(B))$
so that the functional $\bfh(t,\cdot):\bfSp\to\R$ can be continuously extended
on ${\mathcal Y}(\varOmega;B)$, namely
$\nu\mapsto\int_\varOmega\int_Bh(t,x,z)\,\nu_x(\d z)\d t$.
\end{proof}

The convex compactification of $\bfU_{\rm ad}$ from Proposition~\ref{prop-1} is
coarser than the (non-metrizable) Fattorini's construction mentioned in
Sect.\,\ref{sect-int}.

\TTT

\begin{remark}[{\sl Numerical approximation}]\label{rem-numer}\upshape
  The explicit characterization of convex compactifications may suggest
  some approximation strategies. Just as an example
  in the particular case \eqref{special-YM}, one can
  apply the mentioned extreme-point-characterization arguments 
  to see that ${\rm ext}(\calY(I;{\mathcal Y}(\varOmega;B))=
  \{\nu:I\mapsto{\mathcal Y}(\varOmega;B)\ \text{weakly* measurable}\}$
  and then to combine the Kre\u{\i}n-Milman theorem yielding an approximation
  of $\NNU$ in the form $\NNU(t)\sim\sum_i\alpha_i\nu_i(t)$ with
  some $\nu_i:I\to{\mathcal Y}(\varOmega;B)$ and Remark~\ref{rem-approx}
  yielding an approximation
  $\nu_i(t)\sim\sum_j\beta_{ij}\Frakd_{u_{ij}(t)}^{}$ with some
  $u_{ij}(t)\in\bfSp$. The non-negative coefficients satisfies
  $\sum_i\alpha_i=1$ and $\sum_j\beta_{ij}=1$. Altogether,
  $$
  \NNU(t)\sim\sum_i\alpha_i\Big(\sum_j\alpha_i\beta_{ij}\Frakd_{u_{ij}(t)}^{}\Big)
  =\sum_{ij}a_{ij}\Frakd_{u_{ij}(t)}^{}\ \ \text{ with }\ \ a_{ij}=\alpha_i\beta_{ij}\,.
  $$ 
Note that $\sum_{ij}a_{ij}=\sum_i\alpha_i(\sum_j\beta_{ij})=\sum_i\alpha_i=1$
and we obtain a mix of controls $[u_{ij}(t)](x)\equiv u_{ij}(t,x)$
in the spirit of Remark~\ref{rem-approx}. This approximation is convex,
the analytical details about such approximation deserving still some
investigation. Of course, one can also think about combination of
some interpolation over time of some convex combination with
coefficients depending on $t\in I$.
This expectedly opens wide menagerie of possible numerical strategies,
  which remains out of the scope of this article, however.
\end{remark}

\EEE

\begin{remark}[{\sl One generalization}]\label{rem-DM}\upshape
The above construction can be generalized for $B\subset\R^m$ unbounded
by considering a general $\bfSp$ bounded in $L^p(\varOmega;\R^m)$ with some
\TTT specific \EEE $1\le p<\infty$ fixed but not necessarily bounded in
$L^\infty(\varOmega;\R^m)$.
Instead of the conventional Young measures ${\mathcal Y}(\varOmega;B)$, we can
then consider so-called DiPerna-Majda measures
$\TTT\mathrm{DM}\EEE_{\mathcal R}^p(\varOmega;\R^m)$ induced by test-functions of
the form
$g(x)v(z)(1{+}|z|^p)$ with $g\in C(\barOmega)$ and $v$ ranging over some
complete separable ring $\clR\subset C(\R^m)$ containing constants.
\TTT More specifically, $\mathrm{DM}_{\mathcal R}^p(\varOmega;\R^m)$ is
a convex subset of the Radon measures on
$\barOmega\times\gamma_{\mathcal R}\R^m$ attainable from $\bfSp$ when embedded
into $\barOmega\times\gamma_{\mathcal R}\R^m$ via the mapping
$\delta:u\mapsto(h\mapsto\int_\varOmega h(x,u(x))(1{+}|u(x)|^p)\,\d x)$ with
$h\in C(\varOmega)\otimes{\mathcal R}$. \EEE The embedding
$\delta:L^p(\varOmega;\R^m)\to\TTT\mathrm{DM}\EEE_{\clR}^p(\varOmega;\R^m)$ is
homeomorphical and $\TTT\mathrm{DM}\EEE_{\clR}^p(\varOmega;\R^m)$ is convex,
metrizable, and locally compact, having all extreme points of the form
of Diracs $\DELTA_{\bfs}$ with $\bfs\in\bfSp$, cf.\
\cite{kruzik-tr-1,KruRou99SGPS} for $B=\R^m$. As $\bfSp$ is bounded in
$L^p(\varOmega;\R^m)$, the closure of $\DELTA(\bfSp)$ is compact. Like
${\mathcal Y}(\varOmega;B)$, we thus obtained a convex, metrizable separable
compact with extreme points being Diracs, so that the above arguments can be
adopted to this situation, too.
\end{remark}

\section{Application: optimal control of parabolic systems}\label{sec-appl}

Let us briefly outline application to an optimal control of a system of $n$
semilinear parabolic
differential equations. We confine ourselves on homogeneous Dirichlet
conditions, and use the standard notation $H_0^1(\varOmega;\R^n)$ for the
Sobolev space of functions $\varOmega\to\R^n$ whose distributional derivative
is in $L^2(\varOmega;\R^{d\times n})$ and traces on the boundary $\varGamma$
of $\varOmega$ are zero, and similarly $\HH^1(I;X)$ is a Bochner-Sobolev space of
functions $I\to X$ whose distributional derivative is in the Bochner space
$\LL^2(I;X)$. The dual of $H_0^1(\varOmega;\R^n)$ is denoted standardly as
$H^{-1}(\varOmega;\R^n)$. Moreover, we will use the abbreviation 
$\HH^1(I;V,V^*)=\LL^2(I;V)\,\cap\,\HH^1(I;V^*)$ for a Banach space $V$. We will
use it for $V=H_0^1(\varOmega;\R^n)$. We will also use the notation
${\mathcal L}(V)$ for the space of linear bounded operators
from $V$ to $V\cong V^*$. We then consider an initial-boundary value problem
(with $\varGamma$ denoting the boundary of $\varOmega$ with the unit
outward normal $\vec{n}$):
\begin{align}
\!\!\left.\!
\begin{array}{ll}
\mbox{Minimize}\!\!\! & 
\displaystyle{\int_0^T\!\!\!\int_\varOmega
\varphi(t,\bfy(t,\cdot),\bfu(t,\cdot))
  \,\d x
  \d t+
\int_\varOmega{}}\phi(\bfy(T))\,\d x \hfill 
\mbox{(cost functional)}
\\[3mm]
\mbox{subject to}\!\!\! & 
\displaystyle{\frac{\partial\bfy}{\partial t}}-{\rm div}
(A\nabla\bfy)
=\bff(t,\bfy(t,\cdot),\bfu(t,\cdot))
\ \mbox{ in }I{\times}\varOmega
\,,\quad \mbox{(state equation)} 
\\[0mm]
&\bfy\,=\,0\hspace*{4.6em}
\mbox{ on }I{\times}\varGamma\,,\hfill \mbox{(boundary condition)}
\\[.1em]
\ & \bfy(0,\cdot)=\bfy_0\ \ \ \ \ \ \ 
\mbox{ on }\varOmega,\hfill \mbox{(initial condition)}
\\[.1em]
\ & [\bfu(t)](x)\in B
\ \ \ \ \ 
\mbox{ for }\ (t,x)\in I\!\times\!\varOmega\,,\ \ 
\hfill\mbox{(control constraints)} 
\\[.1em]
\ & \bfy\in 
\HH^{1}(I;H_0^1(\varOmega;\R^n),H^{-1}(\varOmega;\R^n))
\,,
\ \ \ \bfu\in\LL^p(I;L^p(\varOmega;\R^m))\,\hspace*{-2cm}
\end{array}\!\!
\right\}\hspace*{-.2cm}
\label{PPAR}
\end{align}
with $\varphi:I\times L^2(\varOmega;\R^n)\times L^\infty(\varOmega;\R^m)\to\R$
and $\bff:I\times L^2(\varOmega;\R^n)\times L^\infty(\varOmega;\R^m)\to
L^2(\varOmega;\R^n)$.
In view of Remark~\ref{rem-1}, the spatially nonlocal right-hand side of the
controlled system can involve an integral over $\varOmega$ so that we could
speak rather about a parabolic integro-differential system.

The relaxation by means of the Young-type measures from
Section~\ref{sec-YM-time} (similarly as from \cite{fattorini,fattorini-4})
records fast oscillations in time but not in space, in contrast to conventional
the conventional Young measures on $I{\times}\varOmega$ which record
fast oscillations simultaneously in time and in space.
Also, the former relaxation allows for a bit more comprehensive
optimality conditions than conventional Young measures on $I{\times}\varOmega$,
cf.\ 
Sect.\ \ref{sec-5} below. To perform our
relaxation, we consider a separable sub-ring $\calR$ of $C(\bfSp)$
with $\bfSp$ from 
\eqref{S} as in Sect.\,\ref{sec-YM-time} and 
qualify the nonlinearities involving the control variable as
\begin{subequations}\label{ass}\begin{align}\nonumber
&\forall\bfy\in C(I;L^2(\varOmega;\R^m)),\ \ \bfv\in H^1(\varOmega;\R^m):\ \ \\ 
&\qquad\varphi{\circ}\bfy:(t,\bfs)\mapsto\varphi(t,\bfy(t),\bfs)\in
\LL^1(I;\calR)\ \ \text{ and}
\\&\qquad
\big\langle\bff{\circ}\bfy,\bfv\big\rangle:
(t,\bfs)\mapsto\big\langle\bff(t,\bfy(t),\bfs),\bfv\big\rangle\in
\LL^1(I;\calR)\,.
\label{ass2}\end{align}\end{subequations}
Then \eqref{PPAR} allows for a continuous extension on the set of
the relaxed controls $\calY(I;\gamma_\calR^{}\bfSp)$ from
Proposition~\ref{prop-1} as:
\begin{align}
\left.
\begin{array}{ll}
\mbox{Minimize}\!\!\! & \displaystyle{\int_0^T\!\!\!
  \int_{\gamma_\calR^{^{}}\bfSp}\!\!\!\!
  \overline\varphi(t,\bfy(t,\cdot),\bfs)\NNU_t(\d\bfs)\d t+
\int_\varOmega\phi(\bfy(T))\,\d x}
\\[1.em]
\mbox{subject to} &\forall\bfv\in \HH^1(I;L^2(\varOmega;\R^m))\,\cap\,
\LL^2(I;H_0^1(\varOmega;\R^m)),\ \bfv(T)=0:
\\[-.0em]
&\ \displaystyle{\int_0^T\!\!\!\bigg(\int_\varOmega\!\!
A\nabla\bfy{:}\nabla\bfv-\bfy{\cdot}\frac{\partial\bfv}{\partial t}}\,\d x
\\[-.3em]
&\qquad\quad
\displaystyle{
-\int_{\gamma_\calR^{^{}}\bfSp}\!\!\!\!\!\!\!\!\!
\big\langle\overlinebff(t,\bfy(t),\bfs),\bfv(t)\big\rangle\NNU_t(\d\bfs)
\bigg)\d t=\int_\varOmega\!\bfy_0{\cdot}\bfv(0)\,\d x}\,,
\\[0.4em]
\ & \bfy\in\HH^1(I;H_0^1(\varOmega;\R^n),H^{-1}(\varOmega;\R^n))\,,
\ \ \ \NNU\in\calY(I;\gamma_\calR^{}\bfSp)\,,
\end{array}
\right\}\hspace*{-.2cm}
\label{RPAR}
\end{align}
where $\overline\varphi(t,\bfy(t,\cdot),\cdot)$ and
$\overlinebff(t,\bfy(t),\cdot)$ denote the (uniquely defined) continuous
extension of  $\varphi(t,\bfy(t,\cdot),\cdot)$ and $\bff(t,\bfy(t),\cdot)$,
respectively. The integral identity in \eqref{RPAR} is a weak formulation of
the initial-boundary value in \eqref{PPAR} arisen by applying once
Green formula in space with using the boundary conditions and
by-part integration in time with using the initial condition.

We will further assume the following ``semi-monotonicity'' condition
for $-\bff(t,\cdot,\bfs)$:
\begin{align}\nonumber
\exists a_1\in L^2(I)&\ \ \forall t\in I\ \ \forall\bfr_1,\bfr_2\in H_0^1(\varOmega;\R^n)\ \ \forall\bfs\in L^\infty(\varOmega;\R^m):
\ \ \ \ 
\\&
\int_\varOmega\!\big(\bff(t,\bfr_1,\bfs)-\bff(t,\bfr_2,\bfs)\big)
\cdot(\bfr_1{-}\bfr_2)\,\d x
\le a_1(t)\|\bfr_1{-}\bfr_2\|_{L^2(\varOmega;\R^n)}^2\,.
\label{ass-monotone}\end{align}
The metrizability and separability of $\calY(I;\gamma_\calR^{}\bfSp)$
allows for stating well-posedness of the relaxed scheme \eqref{RPAR}
conventionally in terms of sequences:

\begin{proposition}[Well-posedness and correctness of
\eqref{RPAR}]\label{prop-2}
  Let  \eqref{ass}, 
\eqref{ass-monotone}, $\phi\in C(L^2(\varOmega;\R^n))$,
 $A\in\R^{(n\times n)^2}$ be positive definite,
and $\bfy_0\in L^2(\varOmega;\R^n)$. Then:
\begin{enumerate}
\item
\eqref{RPAR} possesses a solution and $\min\eqref{RPAR}=\inf\eqref{PPAR}$.
\item
Any infimizing sequence $\{\bfu_k\}_{k\in\N}$ for \eqref{PPAR}
contains a subsequence which, when embedded into
$\calY(I;\gamma_\calR^{}\bfSp)$ by $\DELTA$, converges to some
$\NNU$. Any such limit $\NNU$ solves the relaxed problem \eqref{RPAR}.
\item
Any solution $\NNU$ to \eqref{RPAR} is attainable
by an infimizing sequence $\{\bfu_k\}_{k\in\N}$ for \eqref{PPAR}
in the sense $\NNU=\text{\rm w*-}\lim_{k\to\infty}\DELTA(\bfu_k)$.
\end{enumerate}
\end{proposition}

\begin{proof}[Sketch of the proof]
From positive definiteness of $A$, \eqref{ass}, and
$\bfy_0\in L^2(\varOmega;\R^n)$, we get existence of weak solution $\bfy$
of the initial-boundary value in \eqref{PPAR}. Note that \eqref{ass2}
ensures that all integrals in the integral identity in
\eqref{RPAR} have a good sense. From \eqref{ass-monotone},
we get also uniqueness of this response. 

 This unique solution thus determines a control-to-state mapping
 $\bfpi:\bfu\mapsto\bfy$ from $\LL^p(I;L^p(\varOmega;\R^m))$ to
 $\HH^1(I;H_0^1(\varOmega;\R^n),H^{-1}(\varOmega;\R^n))$. 
Thanks to \eqref{ass2}, this mapping admits a (weak*,weak)-continuous
extension $\overline{\bfpi}:\NNU\mapsto\bfy$ from
$\calY(I;\gamma_\calR^{}\bfSp)$ to
$\HH^1(I;H_0^1(\varOmega;\R^n),H^{-1}(\varOmega;\R^n))$ with
$\bfy$ being the unique weak solution from the integral identity in
\eqref{RPAR}. By the positive definiteness of $A$, the mapping
$\overline{\bfpi}$ is also (weak*,strong)-continuous from
$\calY(I;\gamma_\calR^{}\bfSp)$ to $\LL^2(I;H_0^1(\varOmega;\R^n))$.

We can then view the problems \eqref{PPAR} and \eqref{RPAR} as minimization
problems in terms of the controls only, involving composed functionals
\begin{subequations}\begin{align}\nonumber\\[-1.8em]
&\bfu\mapsto\int_0^T\!\!\int_\varOmega
\varphi(t,[\bfpi(\bfu)](t),\bfu(t))  \,\d x \d t+
  \int_\varOmega{}\!\phi(\bfy(T))\,\d x
  \\[-2.em]\nonumber
\intertext{and its continuous extension}\nonumber\\[-2.em]
&\NNU\mapsto\int_0^T\!\!\int_{\gamma_\calR^{}\bfSp}\!\!\!\!\!
\overline\varphi(t,[\overline{\bfpi}(\NNU)](t),\bfs)\NNU_t(\d\bfs)  \d t+
\int_\varOmega\!\phi(\bfy(T))\,\d x\,,
\label{barJ}\end{align}\end{subequations}
respectively. By density of $\bfU_{\rm ad}$ from \eqref{Uad} in
$\calY(I;\gamma_\calR^{}\bfSp)$ (cf.\ Proposition~\ref{prop-1}) and metrizability
of $\calY(I;\gamma_\calR^{}\bfSp)$, all the assertions 1.--3.\ follow.
\end{proof}

The convexity of $\calY(I;\gamma_\calR^{}\bfSp)$ allows for derivation of
optimality conditions essentially by standard methods of smooth/convex
analysis. This convex geometry directly determines
the resulting so-called maximum principle. For using the standard
smooth analysis and adjoint-equation technique for evaluation of the
G\^ateaux derivative of the composed functional \eqref{barJ},
we assume that, for any $\bfr,\widetilde\bfr\in L^2(\varOmega;\R^n)$,
$\bfs\in\bfSp$, and $t\in I$, it holds
\begin{subequations}\label{ass+}\begin{align}\nonumber
&\varphi(t,\cdot,\bfs):L^2(\varOmega;\R^n)\to\R
    \ \text{ is G\^ateaux differentiable, }\
    \\[-.2em]&\nonumber\qquad
    \forall\,\widetilde\bfy\in\LL^2(I;L^2(\varOmega;\R^n)){:}\ \ 
    \langle\varphi_{\bfr}'{\circ}\bfy,\widetilde\bfy\rangle:(t,\bfs)\mapsto
\langle\varphi_{\bfr}'(t,\bfy(t),\bfs),\widetilde\bfy(t)\rangle\in \LL^1(I;\calR),
\\&\nonumber\qquad\|\varphi_{\bfr}'(t,\bfr,\bfs)\|_{L^2(\varOmega;\R^n)}^{}\le
a_1(t)\|\bfr\|_{L^2(\varOmega;\R^n)}^{}\,,\ \text{  and}
\\&\qquad\|\varphi_{\bfr}'(t,\bfr,\bfs)-\varphi_{\bfr}'(t,\widetilde\bfr,\bfs)\|
_{L^2(\varOmega;\R^n)}^{}\le a_1(t)\|\bfr{-}\widetilde\bfr\|_{L^2(\varOmega;\R^n)}^{}\,,
\\
&\phi:L^2(\varOmega;\R^n)\to\R\ \text{ is G\^ateaux differentiable,}
\\\nonumber
&\bff(t,\cdot,\bfs):L^2(\varOmega;\R^n)\to L^2(\varOmega;\R^n)
\ \text{ is G\^ateaux differentiable,}
   \\&\nonumber\qquad
    \forall\,\widetilde\bfz\in\LL^2(I;{\mathcal L}(L^2(\varOmega;\R^n))){:}\ \ 
    \langle\bff_{\bfr}'{\circ}\bfy,\widetilde\bfz\rangle:(t,\bfs)\mapsto
 \langle\bff_{\bfr}'(t,\bfy(t),\bfs),\widetilde\bfz(t)\rangle\in \LL^1(I;\calR),
\\&\nonumber\qquad\|\bff_{\bfr}'(t,\bfr,\bfs)\|_{{\mathcal L}(L^2(\varOmega;\R^n))}^{}\le
a_2(t)\|\bfr\|_{L^2(\varOmega;\R^n)}^{}\,,\ \text{  and}
\\&\qquad\|\bff_{\bfr}'(t,\bfr,\bfs)-\bff_{\bfr}'(t,\widetilde\bfr,\bfs)\|
_{{\mathcal L}(L^2(\varOmega;\R^n))}^{}\le a_2(t)\|\bfr{-}\widetilde\bfr\|_{L^2(\varOmega;\R^n)}^{}\,,\label{ass+c}
\end{align}\end{subequations}
with $a_1\in L^1(I)$ and $a_2\in L^2(I)$; actually,
a bit more general assumptions would work, too, cf.\
\cite[Sect.\,4.5]{Roub20ROTV}.

\begin{proposition}[Maximum principle for \eqref{RPAR}]\label{prop-4}
  Let \eqref{ass} and \eqref{ass+} hold. Then, any solution
  $\NNU\in\calY(I;\gamma_\calR^{}\bfSp)$ to \eqref{RPAR} satisfies
\begin{align}\nonumber
&\int_{\gamma_\calR^{^{}}\bfSp}\!\!\!\!h_{\bfy,\CHI}(t,\bfs)\,\NNU_t(\d\bfs)
  =\sup_{\bfs\in\bfSp}h_{\bfy,\CHI}(t,\bfs)\ \ \ \text{ for a.a. }\ t\in I
\\\label{max-princ}
&\qquad\qquad\qquad\qquad\text{ with }\ h_{\bfy,\CHI}(t,\bfs)
=\langle\bff(t,\bfy(t),\bfs),\CHI(t)\rangle-\varphi(t,\bfr,\bfs)\,,
\end{align}
with $\bfy=\overline{\bfpi}(\NNU)$ and with
$\CHI\in \HH^{1}(I;H_0^1(\varOmega;\R^n),H^{-1}(\varOmega;\R^n))$ being a weak
solution to the adjoint terminal-boundary-value parabolic problem
\begin{subequations}\label{adj-eq}\begin{align}
&\frac{\partial\CHI}{\partial t}+{\rm div}(A^\top\nabla\CHI)
+\int_{\gamma_\calR^{^{}}\bfSp}\!\!\!\!
\big[h_{\bfy,\CHI}\big]_{\bfr}'(t,\bfs)^\top\,\NNU_t(\d\bfs)=0
&&
\mbox{in }I{\times}\varOmega\,,&&
\\[-.7em]
& 
\CHI\,=\,0 
&&\mbox{on }I{\times}\varGamma\,,
\\&\CHI(T)=\phi_{\bfr}'(\bfy(T))
&&\mbox{on }\varOmega\,.
\end{align}\end{subequations}
\end{proposition}

\begin{proof}[Sketch of the proof] Let us define the extensions
$\overline{\varphi}:I\times H_0^1(\varOmega;\R^n)\times\rca(\gamma_\calR^{}\bfSp)\to
\R\cup\{+\infty\}$
and $\overlinebff:I\times H_0^1(\varOmega;\R^n)\times\rca(\gamma_\calR^{}\bfSp)\to
H^{-1}(\varOmega;\R^n)$ of $\varphi$ and 
$\bff$ by
\begin{subequations}\label{extension}\begin{align}
&\overline{\varphi}(t,\bfr,\NNU)=\begin{cases}
\int_{\gamma_\calR^{}\bfSp}\varphi(t,\bfr,\bfs)\NNU(\d\bfs)&\text{if }\ \NNU\in\prca(\gamma_\calR^{}\bfSp),\\
\qquad+\infty&\text{if }\ \NNU\in\rca(\gamma_\calR^{}\bfSp)\setminus\prca(\gamma_\calR^{}\bfSp)\,,
\end{cases}
\\&\big\langle\overlinebff(t,\bfr,\NNU),\bfv\big\rangle=
\!\int_{\gamma_\calR^{^{}}\bfSp}\!\!\!\!
\big\langle\bff(t,\bfr,\bfs),\bfv\big\rangle\,\NNU(\d\bfs)
\end{align}\end{subequations}
for any $\bfv\in H_0^1(\varOmega;\R^n)$, respectively. By the assumptions
(\ref{ass+}a,b), the functional
$\bfy\mapsto\int_0^T\!\!\int_{\gamma_\calR^{}\bfSp}\overline\varphi(t,\bfy(t),\bfs)\NNU_t(\d\bfs)\,\d t
+\phi(\bfy(T))$ on $\HH^1(I;H_0^1(\varOmega;\R^n),H^{-1}\varOmega;\R^n))$
is G\^ateaux differentiable. Similarly, (\ref{ass+c}) gives smoothness
(namely continuous G\^ateaux differentiability) of
$\bfy\mapsto\overlinebff(\bfy,\NNU)$. Let us further define the extended
composed cost functional
$J:\LL^\infty_{\rm w*}(I;\rca(\gamma_\calR^{}\bfSp))\to\R\cup\{+\infty\}$ defined by
$$
J(\NNU)=\int_0^T\!\!\overline{\varphi}(t,\bfy_{\NNU}(t),\NNU_t)\,\d t
+\phi(\bfy_{\NNU}(T))
$$
with $\bfy_{\NNU}$ being the solution to the controlled
system $\frac{\d\bfy}{\d t}+{\cal A}\bfy=\overlinebff(\bfr,\NNU)$ with
${\cal A}=-{\rm div}(A\nabla\bfy)$ and with $\bfy(0)=\bfy_0$. The functional $J$
has a smooth part determined by $\varphi$, $\bff$, $\bfb$, and $\phi$, and
a nonsmooth but convex part as an indicator function of the convex
subset $\calY(I;\gamma_\calR^{}\bfSp)$ of
$\LL^\infty_{\rm w*}(I;\rca(\gamma_\calR^{}\bfSp))$. The subdifferential $\partial J$
of $J$ can be calculated by the adjoint-equation techniques, leading to
$\partial J(\NNU)=N_{\calY(I;\gamma_\calR^{}\bfSp)}(\NNU)-h_{\bfy_{\NNU},\CHI}$ with
$h_{\bfy_{\NNU},\CHI}\in \LL^1(I;\calR)$ from \eqref{max-princ} and with
$N_{\calY(I;\gamma_\calR^{}\bfSp)}$ denoting the normal cone to
$\prca(\gamma_\calR^{}\bfSp)$
and with $\CHI$ satisfying the integral identity
\begin{align}\nonumber
\!\!\!\int_0^T\!\!\bigg(\big\langle A^\top\nabla\CHI(t),\nabla\bfv(t)\big\rangle
+\Big\langle\CHI(t),\frac{\d\bfv}{\d t}\Big\rangle
+\int_{\gamma_\calR^{^{}}\bfSp}\!\!\!\!
\big\langle[\bff{\circ}\bfy]_{\bfr}'(t,\bfs)^\top\CHI(t),
\bfv(t)\big\rangle\,\NNU_t(\d\bfs)
\bigg)\d t\!\!
\\[-.3em]\label{adj-eq-weak}
=\int_0^T\!\!\!\int_{\gamma_\calR^{^{}}\bfSp}\!\!\!\!\big\langle[\varphi{\circ}\bfy]_{\bfr}'(t,\bfs),\bfv(t)\big\rangle\,\NNU_t(\d\bfs)\d t
+\big\langle\phi_{\bfr}'(\bfy(T)),\bfv(T)\big\rangle\,.\!\!
\end{align}
This is the weak formulation of the terminal-boundary-value problem
\eqref{adj-eq}. The optimality condition $\partial J(\NNU)\ni0$
reads as $\langle\widetilde\NNU{-}\NNU,h_{\bfy_{\NNU},\CHI}\rangle\le0$
for any $\widetilde\NNU\in\calY(I;\gamma_\calR^{}\bfSp)$, i.e.\
$\langle\NNU,h_{\bfy_{\NNU},\CHI}\rangle=\max_{\widetilde\NNU\in\calY(I;\gamma_\calR^{}\bfSp)}
\langle\widetilde\NNU,h_{\bfy_{\NNU},\CHI}\rangle$.
By the density of $\DELTA(\bfU_{\rm ad})$ in
$\calY(I;\gamma_\calR^{}\bfSp)$, this condition
just gives \eqref{max-princ}.
\end{proof}

Exploiting the maximum principle, one can weaken the convexity condition
\eqref{orientor} by considering a smaller set than $\bfSp$ in
\eqref{(20)} excluding arguments which surely cannot satisfy
the maximum principle, cf.\ 
\cite{MunPed01RERN} where the relaxed problems were exploited
but for optimal control of ordinary differential equations. Thus
existence for \eqref{PPAR} can be proved even for nonconvex orientor
fields, cf.\ \cite{RouSch01EOCP}.

\begin{remark}[{\sl Constancy of the Hamiltonian along optimal
      trajectories}]\label{rem-Ham-const}
\upshape
Still one mo\-re condition is sometimes completing the maximum
principle for evolution systems, namely that the
Hamiltonian is constant in time. Here, it is expected that the
augmented Hamiltonian $$
h_{\bfy,\CHI}^{\cal A}(t,\NNU):=
\langle\overlinebff(t,\bfy(t),\NNU)-{\cal A}\bfy(t),\CHI(t)\rangle-
\overline{\varphi}(t,\bfy(t),\NNU)
$$
is constant in
time for any optimal pair $(\bfy,\NNU)$ with $\CHI$ solving \eqref{adj-eq},
i.e.\ the function
$t\mapsto\int_{\bfSp}h_{\bfy,\CHI}^{\cal A}(t,\bfs)\NNU_t(\d\bfs)$ is constant on $I$.
This actually holds only for autonomous systems,
i.e.\ $\varphi$, $\bff$, and $\bfb$ independent of time.
Then, by the following (formal) calculations (with the $t$-variable
not explicitly written), we have
\begin{align}\nonumber
\!\frac{\d}{\d t}h_{\bfy,\CHI}^{\cal A}(t,\NNU)
&=\Big\langle\overlinebff(\bfy,\NNU)-{\cal A}\bfy,\frac{\d\CHI\!}{\d t}\Big\rangle
+\big\langle\overlinebff_t'(\bfy,\NNU),\CHI\big\rangle
-\overline{\varphi}_t'(\bfy,\NNU)
-\Big\langle\overline{\varphi}_{\bfr}'(\bfy,\NNU),\frac{\d\bfy}{\d t}\Big\rangle
\\[-.1em]&\qquad\nonumber
+\Big\langle\big(\overlinebff_{\bfr}'(\bfy,\NNU)-{\cal A}\big)\frac{\d\bfy}{\d t}
,\CHI\Big\rangle
+\Big\langle\frac{\d\NNU}{\d t},
\!\!\!\!\!\lineunder{
h_{\CHI}\circ\bfy
-N_{\prca(\gamma_\calR^{}\bfSp)}^{}(\NNU)_{_{_{_{_{_{}}}}}}\!\!\!\!}{=0 by \eqref{max-princ}}\!\!\!
\Big\rangle
\\[-.8em]&
=\big\langle\overlinebff_t'(\bfy,\NNU),\CHI\big\rangle
-\overline{\varphi}_t'(\bfy,\NNU)\,,
\label{max-Ham-calc}\end{align}
where we
used that $\langle\frac{\d}{\d t}\NNU,N_{\prca(\gamma_\calR^{}\bfSp)}^{}(\NNU)\rangle
=0$. Also we used 
$\frac{\d}{\d t}\bfy+{\cal A}\bfy=\overlinebff(\bfy,\NNU)$
and the adjoint equation (\ref{adj-eq}a,b) in the form
$\frac{\d}{\d t}\CHI-{\cal A}^*\CHI=\overline{\varphi}_{\bfr}'(\bfy,\NNU)
-\overlinebff_{\bfr}'(\bfy,\NNU)^*\CHI$,
which yields
\begin{align}\nonumber
&\Big\langle\overlinebff(\bfy,\NNU)-{\cal A}\bfy,\frac{\d\CHI\!}{\d t}\Big\rangle
+\Big\langle\big(\overlinebff_{\bfr}'(\bfy,\NNU)-{\cal A}\big)\frac{\d\bfy}{\d t}
,\CHI\Big\rangle-\Big\langle\overline{\varphi}_{\bfr}'(\bfy,\NNU),\frac{\d\bfy}{\d t}\Big\rangle
\\&\nonumber=\Big\langle\frac{\d\bfy}{\d t},
{\cal A}^*\CHI{+}\overline{\varphi}_{\bfr}'(\bfy,\NNU)
{-}\overlinebff_{\bfr}'(\bfy,\NNU)^*\CHI\Big\rangle
+\Big\langle\big(\overlinebff_{\bfr}'(\bfy,\NNU){-}{\cal A}\big)\frac{\d\bfy}{\d t}
,\CHI\Big\rangle-\Big\langle\overline{\varphi}_{\bfr}'(\bfy,\NNU),\frac{\d\bfy}{\d t}\Big\rangle=0\,.
\end{align}
From \eqref{max-Ham-calc}, we can see
that $h_{\bfy,\CHI}^{\cal A}$ is constant in time if both $\bff_t'=0$,
$\bfb_t'=0$, and $\varphi_t'=0$. 
\end{remark}

\section{
Some other relaxation schemes}\label{sec-5}

Sometimes, the relaxed problem uses \TTT the \EEE conventional Young measures
from ${\cal Y}(I{\times}\varOmega;B)$. This coarser compactification may
naturally record fast oscillations of infimizing controls both in time and space
simultaneously. 

\def\GG{\widehat\varphi}
\def\FF{f}
\def\PHI{h}

In view of Example~\ref{rem-1}, we can consider rather general
nonlinearities. \TTT To avoid \EEE
too many notational complications,
we consider for example the problem:
\begin{align}
\!\!\left.\!
\begin{array}{ll}
\mbox{Minimize}\!\!\! & 
\displaystyle{\int_0^T\!
  \sum_{i=1}^k\prod_{j=1}^l\GG_{ij}\bigg(t,\!\int_\varOmega\!
\PHI_{ij}(t,x,y(t,x),u(t,x))  \,\d x\!\bigg)
  \d t}
\\[-.5em]&\hspace{15em}+
\displaystyle{\int_\varOmega{}}\!\phi(x,y(T,x))\,\d x 
\\[3mm]
\mbox{subject to}\!\!\! & 
\displaystyle{\frac{\partial y}{\partial t}}-{\rm div}
(A\nabla y)=\displaystyle{
  \FF\big(t,x,y(t,x),u(t,x)  \big)}
\ 
\mbox{ in }I{\times}\varOmega
\,,
\\[0mm]
& y\,=\,0\hspace*{4.6em}
\mbox{ on }I{\times}\varGamma\,,
\\[.2em]
\ & y(0,\cdot)=y_0\ \ \ \ \ \ \ \, 
\mbox{ on }\varOmega\,,
\\[.2em]
\ & u(t,x)\in B
\ \ \ \ \ \ \ \:
\mbox{ for }\ (t,x)\in I\!\times\!\varOmega\,,\ \ 
\\[.2em]
\ & y\in \HH^{1}(I;H_0^1(\varOmega;\R^n),H^{-1}(\varOmega;\R^n))\,,
\ \ \ u\in
L^p(I{\times}\varOmega;\R^m)\,\hspace*{.2cm}
\end{array}\!\!
\right\}\hspace*{-.2cm}
\label{PPAR+}
\end{align}
with $\GG_{ij}:I{\times}\R\to\R$,
$\PHI_{ij}:I{\times}\varOmega{\times}\R^n{\times}
\R^m\to\R$, $\phi:\varOmega{\times}\R^n\to\R$, and
$\FF:I{\times}\varOmega{\times}\R^n{\times}\R^m
   \to\R^n$, $i=1,...,k$ and $j=1,...,l$ with $k,l\in\N$.
This falls into the form \eqref{PPAR} when taking
\begin{subequations}\label{ass++}\begin{align}
&\varphi(t,\bfr,\bfs)
=\sum_{i=1}^k\prod_{j=1}^l\GG_{ij}\bigg(t,\!\int_\varOmega\!\PHI_{ij}(t,x,\bfr(x),\bfs(x))\,\d x\bigg)\,,
\\[-.3em]&\big[\bff(t,\bfr,\bfs)\big](x)=
\FF\big(t,x,\bfr(x),\bfs(x)\big)\,.
\end{align}\end{subequations}
The natural (although not the weakest possible) qualification of these data is
\begin{subequations}\label{ass+++}\begin{align}
 &\PHI_{ij}\in L^1(I{\times}\varOmega;C(\R^n{\times}B))\,,\ \ \ \
 \phi\in L^1(\varOmega;C(\R^n))\,,\ \ \ \
 \\&\GG_{ij}\in \LL^1(I;C(\R))\,,\ \ \text{ and }\ \
 \FF\in L^1(I{\times}\varOmega;C(\R^n{\times}B  )^n) \,.
\end{align}\end{subequations}
Under these assumptions, \eqref{PPAR+} bears
an extension to the conventional Young measures
${\cal Y}(I{\times}\varOmega;B)$, which leads to the relaxed
problem
\begin{align}
\!\!\!\!\left.\!
\begin{array}{ll}
\mbox{Minimize}\!\!\!\!\! & 
\displaystyle{\int_0^T\!\sum_{i=1}^k\prod_{j=1}^l\GG_{ij}\bigg(t,\int_\varOmega\int_B\!
\PHI_{ij}(t,x,y(t,x),z)\,\nu_{t,x}(\d z)  \,\d x\bigg)
  \d t}
\\[-.5em]&\hspace{19em}+
\displaystyle{\int_\varOmega{}}\!\phi(x,y(T,x))\,\d x 
\\[3mm]
\mbox{subject to}\!\!\!\! &
\displaystyle{\int_0^T\!\!\!\!\int_\varOmega\!\!\bigg(
  A\nabla y{:}\nabla v-y{\cdot}\frac{\partial v}{\partial t}
  -\int_B\!\FF\big(t,x,y(t,x),z
\big)\,\nu_{t,x}(\d z)\!\bigg)\,\d x\d t}
\\[.9em]&\hspace*{.7em}\displaystyle{
=\int_\varOmega\!\!y_0{\cdot}v(0)\,\d x
}\qquad\forall v\in \HH^1(I{\times}\varOmega;\R^m)\,,\ \  v(T)=0\,,\ \ 
v|_{I{\times}\varGamma}^{}=0\,,
\\[.8em]
\ & y\in \HH^{1}(I;H_0^1(\varOmega;\R^n),H^{-1}(\varOmega;\R^n))\,,
\ \ \ \nu\in
{\cal Y}(I{\times}\varOmega;B)\,.\hspace*{.2cm}
\end{array}\!\!
\right\}\hspace*{-.2cm}
\label{RPPAR+}
\end{align}
This extension is however not weakly* continuous unless $l=1$ and
all $\GG_{i1}(t,\cdot):\R\to\R$ are affine. 

The resulted (Pontryagin-type) maximum principle is then formulated pointwise
for a.a.\ $(t,x)\in I{\times}\varOmega$. For a very special case $k=1=l$ and
$\GG_{11}(t,\cdot)$ affine, such relaxation scheme has been used
e.g.\ in \cite{Chry06DMSP,ChCoKo??CROM} or also \cite[Sect.4.5.b]{Roub20ROTV}.
In this special case, one can prove also existence of solutions, i.e.\ optimal
relaxed controls from ${\cal Y}(I{\times}\varOmega;B)$.
For a derivation of the mentioned pointwise maximum
principle for the original problem without relaxation in this special
case we refer e.g.\ to 
\cite{casas-1,hu-yong,RayZid99HPPC}. The pointwise constancy of the
Hamiltonian on $I{\times}\varOmega$
however does not seem to hold, in contrast to the finer relaxation
examined before in Remark~\ref{rem-Ham-const}.

In general, the existence of solutions to \eqref{RPPAR+} is however not granted
by usual direct-method arguments unless $l=1$ and $\GG_{i1}(t,\cdot)$ are convex.
In view of \eqref{ass++}, we can exploit also the relaxation scheme
from Section~\ref{sec-appl}. The metrizability and separability of
$\gamma_\calR^{}\bfSp$
allows for a generalization of the (originally finite-dimensional)
Filippov-Roxin \cite{Fili62CQTO,Roxi62ETOC} existence theory for nonconvex
problems. Here we exploit the relaxed problem \eqref{RPAR} similarly as it was
done for finite-dimensional systems in \cite{MunPed01RERN,Roub99CLCE}.

\begin{proposition}[Filippov-Roxin existence for \eqref{RPAR}]\label{prop-3}
  Let the assumptions of Pro\-po\-sition~\ref{prop-2} with $\varphi$ and
  $\bff$ from \eqref{ass++} be fulfilled  and let the so-called
orientor field
\begin{align}\label{orientor}
\bfQ(t,\bfr):=
\big\{(\alpha,\bff(t,\bfr,\bfs))\in \R\times H_0^1(\varOmega;\R^n)^*;\
\alpha\ge\varphi(t,\bfr,\bfs),\ \bfs\in\bfSp\big\}
\end{align}
be convex for a.a.\ $t\in I$ and all $\bfr\in H_0^1(\varOmega;\R^n)$.
Then the following relaxed problem possesses a solution:
\begin{align}
\left.
\begin{array}{ll}
\mbox{Minimize}\!\!\! & \displaystyle{\int_0^T\!\!\!
  \overline\varphi(t,\bfy(t,\cdot),\overlinebfu(t,\cdot))
  \,\d t+\int_\varOmega\phi(\bfy(T))\,\d x}
\\[1.em]
\mbox{subject to} &\forall\bfv\in \HH^1(I;L^2(\varOmega;\R^m))\,\cap\,
\LL^2(I;H_0^1(\varOmega;\R^m)),\ \ \ \bfv(T)=0:
\\[-.0em]
&\ \displaystyle{\int_0^T\!\!\!\bigg(\int_\varOmega\!\!
A\nabla\bfy{:}\nabla\bfv-\bfy{\cdot}\frac{\partial\bfv}{\partial t}}\,\d x
\\[-.7em]
&\qquad\qquad\qquad
\displaystyle{
-
\big\langle\overlinebff(t,\bfy(t),\overlinebfu(t,\cdot)),\bfv(t)\big\rangle
\bigg)\d t=\int_\varOmega\!\bfy_0{\cdot}\bfv(0)\,\d x}\,,
\\[0.2em]
\ & \bfy\in \HH^1(I;H_0^1(\varOmega;\R^n),H^{-1}(\varOmega;\R^n))
\,,
\ \ \
\\[0.2em]
\ & \overlinebfu\in \LL^\infty_{\rm w*}(I;\calR^*),\ \ \
\overlinebfu(t)\in\gamma_\calR^{}\bfSp\ \ \text{ for a.a. }\ t\in I
\,.
\end{array}
\right\}\hspace*{-.2cm}
\label{RPAR++}
\end{align}
\end{proposition}

\begin{proof}[Sketch of the proof]
Let us define
\begin{align}\label{orientor+}
  \overlinebfQ(t,\bfr):=\big\{(\alpha,\overlinebff(t,\bfr,\bfs))\in
  \R\times H^{-1}(\varOmega;\R^n);\
\alpha\ge\overline\varphi(t,\bfr,\bfs),\ \bfs\in\gamma_\calR^{}\bfSp\big\}\,.
\end{align}
For a.a.\ $t\in I$ and all $\bfr\in H_0^1(\varOmega;\R^n)$,
the convexity and closedness of $\overlinebfQ(t,\bfr)$ just means
\begin{align}\label{(20)}\nonumber\\[-2.8em]
&
\overline{\rm co}\,\big[\overline\varphi{\times}\overlinebff\,\big](t,\bfr,\gamma_\calR^{}\bfSp)
\subset
\overlinebfQ(t,\bfr)\,
\end{align}
with ``$\,\overline{\rm co}\,$'' denoting the closed convex full.
By \eqref{(20)}, we get 
\begin{align}
&\int_{\bfSp}\big[\varphi{\times}\bff\,\big](t,\bfy(t),\bfs)\NNU_t(\d\bfs)
\in\overline{\rm co}
\big[\overline\varphi{\times}\overlinebff\,\big](t,\bfy(t),\gamma_\calR^{}\bfSp)
\subset\overlinebfQ(t,\bfy(t))\,.
\label{(24)}\end{align}

Taking a solution $\NNU$ to the relaxed problem \eqref{RPAR}, we put
\begin{align}\nonumber
\BFS(t)=\bigg\{\bfs\in\gamma_\calR^{}\bfSp;\ \ &\overline\varphi(t,\bfy(t),\bfs)\le\!
\int_{\gamma_\calR^{}\bfSp}\!\!\varphi(t,\bfy(t),\bfs)\NNU_t(\d\bfs),
\\[-.5em]
&\qquad\ \ 
\bff(t,\bfy(t),\bfs)=\!\int_{\gamma_\calR^{}\bfSp}\!\overlinebff(t,\bfy(t),\bfs)\NNU_t(\d\bfs)\bigg\}\,,
\label{(25)}\end{align}
Obviously, $\BFS(t)$ is closed for a.a.\ $t\in I$. We further show that
it is also non-empty.
Indeed, by \eqref{(20)}, for any
$(\alpha,\bfq)\in \overlinebfQ(t,\bfy(t))$ there is
$\bfs\in\gamma_\calR^{}\bfSp$ such that
$\alpha\ge\overline\varphi(t,\bfy(t),\bfs)$ and
$\bfq=\overlinebff(t,\bfy(t),\bfs)$. Hence, for the particular choice 
\begin{align}
(\alpha,\bfq)=\big(\alpha(t),\bfq(t)\big):=
  \int_{\gamma_\calR^{}\bfSp}\!\!\big[\overline\varphi{\times}\overlinebff\,\big](t,\bfy(t),\bfs)\NNU_t(\d\bfs),
\label{(26)}\end{align}
the inclusion \eqref{(24)} implies that $\alpha(t)\ge\overline\varphi(t,\bfy(t),\bfs)$
and $\bfq(t)=\overlinebff(t,\bfy(t),\bfs)$
for some $\bfs\in\gamma_\calR^{}\bfSp$, hence $\BFS(t)\ne\emptyset$.

Moreover, the multi-valued mapping $\BFS:I
\tto\gamma_\calR^{}\bfSp$ defined by
\eqref{(25)} is
measurable. Indeed, $\NNU$ weakly* measurable and $\overline\varphi$ and $\overlinebff$
Carath\'eodory mappings imply that $q$ from \eqref{(26)} is measurable.
Furthermore, by \cite[Thm.\,8.2.9]{AubFra90SVA}, the level sets 
$t\mapsto\{\bfs\in\gamma_\calR^{}\bfSp;\ \overline\varphi(t,\bfy(t),\bfs)\le \alpha(t)\}$ and 
$t\mapsto\{\bfs\in\gamma_\calR^{}\bfSp;\ \overlinebff(t,\bfy(t),\bfs)=\bfq(t)\}$ are measurable. By 
\cite[Thm.\,8.2.4]{AubFra90SVA}, the intersection of these level sets, 
which is just $\BFS(t)$, is also a measurable multi-valued mapping.

Then, by 
\cite[Thm.\,8.1.4]{AubFra90SVA}, the multi-valued mapping $\BFS$ possesses a
measurable selection $\overlinebfu(t)\in \BFS(t)$; here separability and metrizability
of $\gamma_\calR^{}\bfSp$ were used. 

In view of \eqref{(25)},
$\overlinebff(t,\bfy(t),\overlinebfu(t))=\bfq(t)=\int_{\gamma_\calR^{}\bfSp}\overlinebff(t,\bfy(t),\bfs)\NNU_t(\d\bfs)
$
so that the pair
$(\overlinebfu,\bfy)$ is admissible for \eqref{RPAR}, and moreover 
\begin{align}\nonumber
\int_0^T\!\!\!\overline\varphi(t,\bfy(t),\overlinebfu(t))\,\d t
&\le\!\int_0^T\!\!\!\alpha(t)\,\d t
=\!\int_0^T\!\!\!\int_{\gamma_\calR^{}\bfSp^{^{}}}\!\!\!\!\!\!\!\!\!\overline\varphi(t,\bfy(t),\bfs)\NNU_t(\d\bfs)\,\d t
=\min\eqref{RPAR}=\inf\eqref{RPAR++}\,.
\end{align}
This $\overlinebfu$ thus solves \eqref{RPAR++}.
\end{proof}

The particular choice \eqref{tests} allows for usage of
Lemma~\ref{lem2}. In this case, Proposition~\ref{prop-3} with Lemma~\ref{lem2}
gives $\overlinebfu:I\to{\cal Y}(\varOmega;B)$ as a solution to
\eqref{RPAR++}. Then, in view of special nonlinearities involved in
\eqref{PPAR+}, we can use Lemma~\ref{lem1}, which leads to a relaxation using
a certain Young measure valued on the original set $\bfSp$ from \eqref{S} as
actually used in \eqref{PPAR}, provided we weaken a bit the measurability of
Young measures. More specifically, we define:
\begin{align}\nonumber
{\boldsymbol w\text{-}}\calY(I;\bfSp):=\big\{
\MMU:I\to\prca(\bfSp);\ &\forall
\bfh\in \LL^1(I;C({\cal Y}(\varOmega;B))|_{\bfSp}):
\\[-.2em]&\,t\mapsto\big\langle\MMU_t,\bfh(t,\cdot)\big\rangle
  \text{ is measurable}\big\}\,,
  \label{w-Y}\end{align}
where we used again the convention $\MMU_t:=\MMU(t)$, so that we will write
$\MMU=\{\MMU_t\}_{t\in I}$ in what follows. We call elements of
${\boldsymbol w\text{-}}\calY(I;\bfSp)$ as weak-Young measures. Note that the
set of test functions in \eqref{w-Y} is smaller than the nonseparable
space $\LL^1(I;C(\bfSp))$ and thus weak-Young measures do not live in
$\calY(I;\bfSp)$ in general.

\begin{corollary}
Let the assumptions of Pro\-po\-sition~\ref{prop-2} with $\calR$ from
\eqref{tests} and \eqref{orientor} hold
for $\varphi$ and $\bff$ from \eqref{ass++}. Then there exists a solution to
the following relaxed problem:
\begin{align}
\left.
\begin{array}{ll}
\mbox{Minimize}\!\!\! & \displaystyle{\int_0^T\!
  \sum_{i=1}^n\prod_{j=1}^m\GG_{ij}\bigg(t,\! \int_{\bfSp^{^{}}}\!\!\!
\big[\varPsi([\PHI{\circ}\bfy](t))\big](u)\,\MMU_t(\d u) \bigg)\d t+
\int_\varOmega\!\phi(\bfy(T))\,\d x}
\\[1.em]
\mbox{subject to}\!\! &\forall\bfv\in \HH^1(I;L^2(\varOmega;\R^m))\,\cap\,
\LL^2(I;H_0^1(\varOmega;\R^m)),\ \ \ \bfv(T)=0:
\\[-.0em]
&\ \displaystyle{\int_0^T\!\!\!\bigg(\int_\varOmega\!\!
A\nabla\bfy{:}\nabla\bfv-\bfy{\cdot}\frac{\partial\bfv}{\partial t}}\,\d x
\\[-.4em]
&\qquad\quad
\displaystyle{
-\int_{\bfSp^{^{}}}\!\!\!
\big[\varPsi(\langle [\bff{\circ}\bfy](t),\bfv(t)\rangle)\big]
      (u)\,\MMU_t(\d u)
\bigg)\d t=\int_\varOmega\!\bfy_0{\cdot}\bfv(0)\,\d x}\,,
\\[0.2em]&
\MMU_t\in{\rm srca}_1^+(\bfSp)\ \text{ for a.a.}\ t\in I\,,
\\[0.2em]
\ & \bfy\in 
\HH^1(I;H_0^1(\varOmega;\R^n),H^{-1}(\varOmega;\R^n))
\,,
\ \ \ \MMU\in{\boldsymbol w\text{-}}\calY(I;\bfSp)\,,
\end{array}
\right\}\hspace*{-.2cm}
\label{RPAR+++}
\end{align}
where $\varPsi$ is from \eqref{Psi}. Moreover, if also
\begin{subequations}\label{ass++++}\begin{align}
 &[\PHI_{ij}]_r'\in L^1(I{\times}\varOmega;C(\R^n{\times}B)^n)\,,\ \ \ \
 \phi_r'\in L^1(\varOmega;C(\R^n)^n)\,,\ \ \ \
 \\&[\GG_{ij}]_r'\in \LL^1(I;C(\R))\,,\ \ \text{ and }\ \
 [\FF]_r'\in L^1(I{\times}\varOmega;C(\R^n{\times}B )^{n\times n}) \,,
\end{align}\end{subequations}
then this solution satisfies, for a.a.\ $t\in I$, the maximum principle 
\begin{align}\label{max+}
  \int_{\bfSp^{^{}}}\!\!\!
  h_{\bfy,\CHI}(t,u)\,\MMU_t(\d u)=\sup_{u\in\bfSp}h_{\bfy,\CHI}(t,u)
\end{align}
with $h_{\bfy,\CHI}$ from \eqref{max-princ} with $\CHI$ satisfying
the adjoint terminal-boundary-value parabolic problem, 
written in the weak form here as
\begin{align}\nonumber
&\!\!\!\int_0^T\!\!\bigg(\big\langle A^\top\nabla\CHI(t),\nabla\bfv(t)\big\rangle
+\Big\langle\CHI(t),\frac{\d\bfv}{\d t}\Big\rangle
+\!\int_{\bfSp}\!\!\!
\big[\varPsi(\langle[\bff{\circ}\bfy]_{\bfr}'(t)^\top\CHI(t),
\bfv(t)\rangle)\big](u)\,\MMU_t(\d u)
\bigg)\d t\!\!
\\[-.3em]\label{adj-eq-weak+}
&\hspace*{6em}=\int_0^T\!\!\!\int_{\bfSp}\!\!\!\big[\varPsi(\big\langle[\varphi{\circ}\bfy]_{\bfr}'(t),\bfv(t)\big\rangle)\big](u)\,\MMU_t(\d u)\,\d t
+\big\langle\phi_{\bfr}'(\bfy(T)),\bfv(T)\big\rangle
\end{align}
for all $\bfv\in \HH^1(I;L^2(\varOmega;\R^m))\,\cap\,
\LL^2(I;H_0^1(\varOmega;\R^m))$ with $\bfv(0)=0$.
\end{corollary}

\begin{proof}
  Recall that now, for the choice \eqref{tests},
  $\gamma_{\calR}^{}\bfSp\cong {\cal Y}(\varOmega;B)$, cf.\
  Lemma~\ref{lem2}. Take $\overlinebfu:I\to{\cal Y}(\varOmega;B)$
  a solution to \eqref{RPAR++}. In particular, for any $\bfh\in \LL^1(I;\calR)$,
  the function
  $t\mapsto\langle\overlinebfu(t),\bfh(t)\rangle
  =\int_\varOmega\int_B\bfh(t,x,z)\overlinebfu(t,\d z)\d x=
   \overlinebfh(t,\overlinebfu(t))$
  is measurable (and integrable); here $\overlinebfh(t,\cdot)$ denotes the
  weakly* continuous extension of
  $\bfh(t,\cdot):\bfSp\to\R$ on ${\cal Y}(\varOmega;B)$.

  The functions $[\PHI_{ij}{\circ}\bfy](t):(x,z)\mapsto
  \PHI(t,x,y(t,x),z)$
  and $\langle [\bff{\circ}\bfy](t),\bfv(t)\rangle:(x,z)\mapsto
  \langle f(t,x,y(t,x),z),v(t,x)\rangle$ belong to $ L^1(\varOmega;C(B))$.
  Therefore, by Lemma~\ref{lem1}, for some $\MMU_t\in\prca(\bfSp)$, we have
\begin{subequations}\label{u=mu}\begin{align}\nonumber\\[-2.5em]
&\hspace*{-4em}\GG_{ij}(t,\langle\overlinebfu(t,\cdot),[\PHI_{ij}{\circ}\bfy](t,\cdot)\rangle)
  =g_{ij}(t,\langle\MMU_t,\varPsi([\PHI_{ij}{\circ}\bfy](t))\rangle)
  \ \text{ and}\hspace*{-2em}
\\&\hspace*{-4em}\big\langle\overlinebff(t,\bfy(t),\overlinebfu(t,\cdot)),\bfv(t)\big\rangle
=\big\langle\MMU_t,\varPsi(\langle [\bff{\circ}\bfy](t),\bfv(t)\rangle)\big\rangle\,.
\end{align}\end{subequations}
  Thus $\min$\eqref{RPAR++}$\,\ge\min$\eqref{RPAR+++}.

On the other hand, also $\min\eqref{RPAR++}\le\min\eqref{RPAR+++}$
because, for any $\MMU$ admissible for \eqref{RPAR+++}, there
is some $\overlinebfu:I\to{\cal Y}(\varOmega;B)$ such that
\eqref{u=mu} holds. Thus such $\overlinebfu$ is admissible for
\eqref{RPAR++}, yielding the cost not lower than $\min\eqref{RPAR+++}$.
Here the definition \eqref{ext+++} of ${\rm srca}_1^+(\bfSp)$ has been used.

By \eqref{ass++++}, it can be seen that \eqref{ass+} for $\varphi$ and $\bff$
from \eqref{ass++} is satisfied with $\calR$ from \eqref{tests}.
Then one can use Proposition~\ref{prop-4}. By this way, 
\eqref{max-princ} results to \eqref{max+} while \eqref{adj-eq-weak}
gives \eqref{adj-eq-weak+}.
\end{proof}

\end{sloppypar}
\end{document}